\documentclass[11pt]{article}
\usepackage{amsmath,amssymb,amsfonts}
\usepackage[latin1]{inputenc}

\usepackage{amsfonts}
\usepackage{color}
\usepackage{epsfig}
 \usepackage{dsfont}
\usepackage{amsmath}
\usepackage{hyperref}

\usepackage{graphicx} 

\usepackage{booktabs} 
\usepackage{array} 
\usepackage{paralist} 
\usepackage{verbatim} 
\usepackage{subfig} 
\usepackage{url}
\usepackage{hyperref}
\usepackage{float}

\usepackage{epsf}
\usepackage{psfrag} 
\usepackage{epsfig,float,psfrag}
\hypersetup{pdfborder={0 0 0},colorlinks,urlcolor=blue}  

\numberwithin{equation}{section}
\newtheorem{thm}{Theorem}[section]

\newtheorem{alem}[thm]{Lemma}
\newtheorem{aprop}[thm]{Proposition}

\newtheorem{arem}[thm]{Remark}
\newenvironment{adem}[1][]%
   {\ \\ {\bf Proof #1~: }}%
   {\hfill\mbox{\rule{2 true mm}{3 true mm}}\vskip 2 ex\noindent}
\newenvironment{_adem}[1][]%
   {\ \\ {\bf Proof #1}}%
   {\hfill\mbox{\rule{2 true mm}{3 true mm}}\vskip 2 ex\noindent}
   {\ \\ {\bf Example #1~: }}%
   {\hfill\mbox{\rule{2 true mm}{3 true mm}}\vskip 2 ex\noindent}

\title{Asymptotic error distribution for the Ninomiya-Victoir scheme in
the commutative case}
\author{A. Al Gerbi, B. Jourdain\thanks{Universit\'e Paris-Est, Cermics (ENPC), INRIA, F-77455, Marne-la-Vall\'ee, France
    e-mails: jourdain@cermics.enpc.fr, anis.al-gerbi@cermics.enpc.fr - This research benefited
    from the support of the ``Chaire Risques Financiers'', Fondation du
    Risque.}~ and E. Cl\'ement\thanks{Universit\'e Paris-Est, LAMA (UMR 8050), UPEMLV, UPEC, CNRS, F-77454, Marne-la-Vall\'ee, France,
    e-mail: emmanuelle.clement@u-pem.fr.}}

\begin{document}
\maketitle

In \cite{NV2} we proved strong convergence with order $1$ of the Ninomiya-Victoir scheme $X^{NV}$ with time step $T/N$ to the solution $X$ of the limiting SDE when the Brownian vector fields commute. In this paper, we prove that the normalized error process $N \left(X - X^{NV}\right)$ converges to an affine SDE with source terms involving the Lie brackets between the Brownian vector fields and the drift vector field. This result ensures that the strong convergence rate is actually $1$ when the Brownian vector fields commute, but at least one of them does not commute with the drift vector field. When all the vector fields  commute the limit vanishes. Our result is consistent with the fact that the Ninomiya-Victoir scheme solves the SDE in this case. 

\section{Introduction}

We consider a general $n$-dimensional stochastic differential equation, driven by a $d$-dimensional standard Brownian motion $W = \left(W^1,\ldots,W^d\right)$, of the form
\begin{equation}
\left\{
    \begin{array}{ll}
dX_t = b(X_t) dt + \sum \limits_{j=1}^d \sigma^j(X_t)dW_t^j, ~	 t \in [0,T]\\
X_0 = x,
\end{array}
\right.
\label{EDS_ITO}
\end{equation}
where $x \in \mathbb{R}^n$ is the starting point, $b: \mathbb{R}^n \longrightarrow \mathbb{R}^n$ is the drift coefficient and $\sigma^j: \mathbb{R}^n \longrightarrow \mathbb{R}^n, j \in \left\{1,\ldots,d\right\}$, are the Brownian vector fields.
We are interested in the study of the normalized error process for the Ninomiya-Victoir scheme. To do so we will consider in the whole paper a regular time grid, with time step $h = T/N$, of the time interval $[0,T]$. We introduce some notations to define the Ninomiya-Victoir scheme. Let
\begin{itemize}
\item $\left(t_k = k h \right)_{k \in [\![0;N]\!]}$ be the subdivision of $[0,T]$ with equal time step $h$,
\item $\Delta W^j_s = W^j_s -W^j_{t_k}$, for $s \in \left(t_k,t_{k+1}\right]$ and $j \in \left\{1,\ldots,d\right\}$,
\item $\Delta s = s - t_k$, for $s \in \left(t_k,t_{k+1}\right]$.
\end{itemize}
For $V: \mathbb{R}^n \longrightarrow \mathbb{R}^n$ Lipschitz continuous, $\exp(tV)x_0$ denotes the solution, at time $t \in \mathbb{R}$, of the following ordinary differential equation in $\mathbb{R}^n$  
\begin{equation}
\left\{
    \begin{array}{ll}
 \frac{dx(t)}{dt}  = V\left(x(t)\right) \\
  x(0) = x_0. 
\end{array}
\right.
\label{ODE}
\end{equation}

To deal with the Ninomiya-Victoir scheme, it is more convenient to rewrite the stochastic differential equation \eqref{EDS_ITO} in Stratonovich form. Assuming $\mathcal{C}^{1}$ regularity for the vector fields, the Stratonovich form of \eqref{EDS_ITO} is given by:
\begin{equation}
\left\{
    \begin{array}{ll}
dX_t = \sigma^0(X_t) dt + \sum \limits_{j=1}^d \sigma^j(X_t)\circ dW_t^j \\
X_0 = x,
\end{array}
\right.
\label{EDS_STO}
\end{equation}
where $\sigma^0 = b - \displaystyle \frac{1}{2} \sum \limits_{j=1}^d \partial \sigma^j \sigma^j $ and $\partial \sigma^j$ is the Jacobian matrix of $\sigma^j$ defined as follows
\begin{equation}
\label{Jacobian}
\partial \sigma^j = \left(\left(\partial \sigma^j \right)_{ik}\right)_{i,k \in [\![1;n]\!] } = \left(\partial_{x_k} \sigma^{ij} \right)_{i,k \in [\![1;n]\!] }. 
\end{equation} 
We recall that the Ninomiya-Victoir scheme \cite{NV} is given by: 
\begin{itemize}
\item starting point: $X^{NV,\eta}_{t_0} = x$,
\item for $k \in \left\{0\ldots,N-1\right\}$,
 if $\eta_{k+1} = 1 $:
\begin{equation}
 X^{NV,\eta}_{t_{k+1}} = \exp\left(\frac{h}{2}\sigma^0\right) \exp\left (\Delta W^d_{t_{k+1}}\sigma^d \right) \ldots \exp\left (\Delta W^1_{t_{k+1}}\sigma^1 \right)  \exp\left(\frac{h}{2}\sigma^0\right) X^{NV,\eta}_{t_{k}}, 
\label{case 1 NV}
\end{equation}
and if $\eta_{k+1} = -1 $:
\begin{equation}
X^{NV,\eta}_{t_{k+1}} = \exp\left(\frac{h}{2}\sigma^0\right) \exp\left (\Delta W^1_{t_{k+1}}\sigma^1 \right) \ldots \exp\left (\Delta W^d_{t_{k+1}}\sigma^d \right)  \exp\left(\frac{h}{2}\sigma^0\right)  X^{NV,\eta}_{t_{k}},
\label{case 2 NV}
\end{equation}
\end{itemize}
where $\eta = \left(\eta_k\right)_{k \ge 1}$ is a sequence of independent, identically distributed Rademacher random variables independent of $W$.
In \cite{NV1}, we proved strong convergence with order $1/2$:
\begin{equation}
\forall p \ge 1, \exists C_{NV} \in \mathbb{R}_+^*, \forall N \in \mathbb{N}^*, \mathbb{E}\left[ \underset{0 \leq k\leq N}{\max}\left\|X_{t_k} - X^{NV,\eta}_{t_k}\right\|^{2p} \right] \leq C_{NV} \left(1+ \left\| x\right\|^{2p}\right) h^p.
\end{equation}
In \cite{NV2}, we studied the stable convergence in law of the normalized error defined by $V^N = \sqrt{N}\left(X - X^{NV,\eta}\right)$.
The theory of stable convergence, introduced by R\'enyi \cite{Renyi} was developed by Kurtz-Protter\cite{KP}, Jacod \cite{Jacod} and Jacod-Protter \cite{JP}. The asymptotic distribution of the normalized error for the Euler continuous scheme was established by Kurtz and Protter in \cite{KP}. The asymptotic behavior of the normalized error processes for the Milstein scheme \cite{Milstein}, which is known to exhibit strong convergence with order $1$, was studied by Yan in \cite{Yan}. In both cases, the normalized error converges to the solution of an affine SDE with a source term involving additional randomness given by a Brownian motion independent of the one driving both the SDE and the scheme. 
In \cite{NV2}, we showed the stable convergence in law of $V^N$ to the solution $V$ of the affine SDE with source terms involving the Lie brackets between the Brownian vector fields:
\begin{equation*}
V_t = \sqrt{\frac{T}{2}} \sum \limits_{j=1}^d \sum \limits_{m=1}^{j-1}  \int_0^t \displaystyle \left[\sigma^j,\sigma^m\right]\left(X_s\right) dB^{j,m}_s + \int_0^t \displaystyle \partial b\left(X_s \right)V_s ds + \sum \limits_{j=1}^d \int_0^t \displaystyle \partial \sigma^j\left(X_s \right)V_s dW_s^j,
\end{equation*} 
where $\left[\sigma^j,\sigma^m\right] = \partial \sigma^m \sigma^j - \partial \sigma^j \sigma^m $, for $j,m \in \left\{1,\ldots,d\right\}, m < j$, denotes the Lie bracket between the Brownian vector fields $\sigma^j$ and $\sigma^m$, $\partial b$ is the Jacobian matrix of $b$, defined analogously to \eqref{Jacobian}, and $\left(B_t\right)_{0\leq t \leq T}$ is a standard $\frac{d(d-1)}{2}$-dimensional Brownian motion independent of $W$. The limit vanishes when the Brownian vector fields commute:
\begin{equation}
\label{Comm_condition}
\forall j,m \in \left\{1,\ldots,d\right\}, \left[ \sigma^j,\sigma^m\right] = \partial \sigma^m\sigma^j - \partial \sigma^j \sigma^m = 0\tag{$\mathcal{C}$}.
\end{equation} 

When the Brownian vector fields $\sigma^j$, for $ j\in \left\{1,\ldots,d\right\}$, commute, the order of integration of these fields no longer matters, since Frobenius' theorem ensures (see \cite{Dieudonee} or \cite{Doss}) the commutativity of the associated flows. The sequence $\eta$ is then useless. Therefore, the Ninomiya-Victoir scheme may be written as follows
\begin{itemize}
\item starting point: $X^{NV}_{t_0} = x$,
\item for $k \in \left\{0\ldots,N-1\right\}$,
\begin{equation}
 X^{NV}_{t_{k+1}} = \exp\left(\frac{h}{2}\sigma^0\right) \exp\left (\Delta W^d_{t_{k+1}}\sigma^d \right) \ldots \exp\left (\Delta W^1_{t_{k+1}}\sigma^1 \right)  \exp\left(\frac{h}{2}\sigma^0\right) X^{NV}_{t_{k}}. 
\label{case 1 NV}
\end{equation}
\end{itemize}
Under some regularity assumptions, we proved, in \cite{NV2}, strong convergence with order $1$ of the Ninomiya-Victoir scheme when the commutativity condition \eqref{Comm_condition} holds. More precisely, we showed the following result.
\begin{thm}
\label{SC2}
Assume that 
\begin{itemize}
\item $\forall j \in \left\{1,\ldots,d\right\}, \sigma^j \in \mathcal{C}^{1}\left(\mathbb{R}^n,\mathbb{R}^n\right)$  with bounded first order derivatives,
\item $\sigma^0 \in \mathcal{C}^{2}\left(\mathbb{R}^n,\mathbb{R}^n\right)$ with bounded first order derivatives and polynomially growing second order derivatives,
\item $\sum \limits_{j=1}^d \partial \sigma^j \sigma^j$ is a Lipschitz continuous function,
\end{itemize}
and that the commutativity condition \eqref{Comm_condition} holds. Then
\begin{equation}
\forall p \ge 1, \exists C^{\prime}_{NV} \in \mathbb{R}_+^*, \forall N \in \mathbb{N}^*, \mathbb{E}\left[ \underset{0 \leq k\leq N}{\max}\left\|X_{t_k} - X^{NV}_{t_k}\right\|^{2p} \right] \leq C^{\prime}_{NV} h^{2p}.
\end{equation}
\end{thm}
In the present paper, we assume that the commutativity condition \eqref{Comm_condition} holds and we focus on the convergence in law of the normalized error defined by $U^N = N\left(X - X^{NV}\right)$. This paper is organized as follows. In section 2, we define an adapted interpolation between time grid points and derive its It\^o decomposition. Then, we provide a suitable decomposition of the normalized error $U^N = N\left(X - X^{NV}\right)$ of the form
\begin{equation}
\label{method_AED}
U^{N}_t = Q^N_t + J^N_t + \left( \int_0^t \displaystyle  H^{0,N}_s U^N_s ds + \sum \limits_{j=1}^d \int_0^t \displaystyle H^{j,N}_s U^N_s dW_s^j \right),
\end{equation} 
where $H^{j,N}, j \in  \left\{0,\ldots,d\right\}$ take values in $\mathbb{R}^n \otimes \mathbb{R}^n$, $Q^N$ is a remainder term and $J^N$ a source term, to study its stable convergence in law.
In section 3, we prove the stable convergence in law of $U^N$ to the solution of the following SDE:
\begin{equation}
U_t = \displaystyle \frac{T}{2\sqrt{3}} \sum \limits_{j=1}^d \displaystyle \int_{0}^{t} \displaystyle \left[ \sigma^0, \sigma^j\right]\left(X_s\right) d\tilde{B}^j_s + \int_0^t \displaystyle \partial b\left(X_s \right)U_s ds + \sum \limits_{j=1}^d \int_0^t \displaystyle \partial \sigma^j\left(X_s \right)U_s dW_s^j
\end{equation} 
where $\left(\tilde{B}_t\right)_{0\leq t \leq T}$ is a standard $d$-dimensional Brownian motion independent of $W$. This result ensures that the strong convergence rate is actually $1$ when the Brownian vector fields commute, but at least one of them does not commute with the drift vector field $\sigma^0$. It is not surprising that the limit vanishes when all the vector fields $\sigma^j$, for $ j\in \left\{0,\ldots,d\right\}$, commute, since the Ninomiya-Victoir scheme solves the SDE \eqref{EDS_ITO} in this case. 

\textbf{Notation}\\\\
In the following, we introduce some more notations which will be used throughout this paper. 
\begin{itemize}
\item Let $\hat{\tau}_s$ be the last time discretization before $s \in [0,T] $, ie $\hat{\tau}_s = t_k$ if $s \in \left(t_k,t_{k+1}\right]$, and for $s = t_0 = 0 $, we set $\hat{\tau}_0 = 0$.
\item Let $\check{\tau}_s$ be the first time discretization after $s \in [0,T]$, ie $\check{\tau}_s = t_{k+1}$ if $s \in \left(t_k,t_{k+1}\right]$, and for $s = t_0 = 0 $, we set $\check{\tau}_0 = 0$.
\item For the vector field $\sigma^j$, $j \in \left\{0,\ldots,d\right\}$, $\partial^2 \sigma^j$ denotes the $n \times n \times n$-tensor $\left(\partial^2 \sigma^j\right)^{i,k,l} = \partial^2_{x_lx_k} \sigma^{ij}$.
\item  The tensor product, between a $m \times p \times q-$tensor A and a vector b in $\mathbb{R}^q$ is denoted by $A\odot b$:
\begin{equation}
\left(A\odot b\right)_{i,k} = \sum \limits_{l=1}^q A^{i,k,l} b^l, \forall i \in \left\{1, \ldots,m\right\}, \forall k \in \left\{1, \ldots,p\right\}.
\end{equation}
\item To lighten up the notation, $\left\| . \right\|$ will denote both the Euclidean norm in $\mathbb{R}^n$ and its associated operator norm in $\mathbb{R}^n \otimes \mathbb{R}^n$.
\item For a finite-dimensional normed vector space $\left(\mathcal{S},\left\| . \right\|_{\mathcal{S}}\right)$, $LIP_{loc}^{pgc}\left(\mathcal{S}\right)$ denotes the space of locally Lipschitz with polynomially growing Lipschitz constant functions from $\mathbb{R}^n $ to $\mathcal{S}$:
$$LIP_{loc}^{pgc}\left(\mathcal{S}\right) = \left\{ F : \mathbb{R}^n \longrightarrow \mathcal{S}, \exists c \in \mathbb{R}_+^*, \exists q \in \mathbb{N},\forall x,y \in \mathbb{R}^n, \left\| F(x) - F(y) \right\|_{\mathcal{S}} \leq c \left(1 + \left\| x \right\|^q \vee \left\| y \right\|^q\right) \left\| x - y \right\|    \right\}.$$
\end{itemize}
\begin{arem} ~
\begin{itemize}
\item If $F \in \mathcal{C}^1 \left(\mathbb{R}^n,\mathcal{S}\right)$ with polynomially growing first order derivatives, then $F \in  LIP_{loc}^{pgc}\left(\mathcal{S}\right)$.
\item If $A \in LIP_{loc}^{pgc}\left(\mathbb{R}^m \otimes \mathbb{R}^p \otimes \mathbb{R}^q\right)$, $M \in LIP_{loc}^{pgc}\left(\mathbb{R}^p \otimes \mathbb{R}^q\right)$ and $b \in LIP_{loc}^{pgc}\left(\mathbb{R}^q \right)$, then $A\odot b \in LIP_{loc}^{pgc}\left(\mathbb{R}^m \otimes \mathbb{R}^p \right)$ and $M b \in LIP_{loc}^{pgc}\left(\mathbb{R}^p \right)$.
\end{itemize}
\end{arem}

\section{Adapted interpolation and main result}
\subsection{Main result}

To study the normalized error process $U^N = N\left( X - X^{NV} \right)$, using the framework described in section 2 of \cite{NV2}, we provide the following adapted interpolation between time grid points:
\begin{equation}
\left\{
    \begin{array}{ll}
 X^{NV}_t  = h_{d+1}\left(\frac{\Delta t}{2}, \Delta W_t,\frac{\Delta t}{2}; X^{NV}_{\hat{\tau}_{t}} \right) \\
  X^{NV}_0 = x,
\end{array}
\right.
\end{equation}
where $h_{d+1}: \mathbb{R}^{d+1} \longrightarrow \mathbb{R}^n$ is defined by
\begin{equation}
\label{eq_hd+1}
h_{d+1}\left(t_0,\ldots,t_{d+1};y\right) = \exp\left(t_{d+1}\sigma^0\right)\exp\left(t_{d}\sigma^d\right)\ldots\exp\left(t_{1}\sigma^1\right)\exp\left(t_{0}\sigma^0\right) y, 
\end{equation}
for the initial condition $y \in \mathbb{R}^n$.
The main result of this paper is the following theorem, which gives the stable convergence in law of the normalized error process $U^N$.
\begin{thm}
\label{T_EP_CC}
Assume that
\begin{itemize}
\item $\forall j \in \left\{0,\ldots,d\right\}, \sigma^j \in \mathcal{C}^{2}\left(\mathbb{R}^n,\mathbb{R}^n\right)$  with bounded first order derivatives, $\partial \sigma^j \in LIP_{loc}^{pgc}\left(\mathbb{R}^n \otimes \mathbb{R}^n\right)$ and $\partial^2 \sigma^j  \in LIP_{loc}^{pgc}\left(\mathbb{R}^n \otimes \mathbb{R}^n \otimes \mathbb{R}^n\right)$,
\item $\forall j \in \left\{1,\ldots,d\right\}, \partial \sigma^j \sigma^j $ is a Lipschitz continuous function,
\end{itemize}
and that the commutativity condition \eqref{Comm_condition} holds.
Then:
\begin{equation}
U^N = N\left( X - X^{NV} \right)  \overset{stably}{\underset{N \to +\infty}{\Longrightarrow}} U
\end{equation}
where $U$ is the unique solution of the following affine equation:
\begin{equation}
U_t = \displaystyle \frac{T}{2\sqrt{3}} \sum \limits_{j=1}^d \displaystyle \int_{0}^{t} \displaystyle \left[ \sigma^0, \sigma^j\right]\left(X_s\right) d\tilde{B}^j_s + \int_0^t \displaystyle \partial b\left(X_s \right)U_s ds + \sum \limits_{j=1}^d \int_0^t \displaystyle \partial \sigma^j\left(X_s \right)U_s dW_s^j
\label{EQ}
\end{equation} 
and $\left(\tilde{B}_t\right)_{0\leq t \leq T}$ is a standard $d$-dimensional Brownian motion independent of $W$.
\end{thm}
In order to prove this theorem we will proceed as follows. Firstly, we will derive the dynamics of the interpolated Ninomiya-Victoir scheme  $X^{NV}$ using It\^o's formula. Then we will decompose the dynamics as follows
 \begin{equation}
\label{decomposition_1_XVN}
dX^{NV}_t = b(X^{NV}_t) dt + \sum \limits_{j=1}^d \sigma^j(X^{NV}_t)dW_t^j + d\bar{J}^N_t + d\bar{Q}^N_t
\end{equation}
where $\bar{J}^N$ is a term with strong order $1$ and $\bar{Q}^N$ is a remainder term such that $N\bar{Q}^N$ uniformly converges in probability to 0. Using \eqref{decomposition_1_XVN} we can derive a convenient decomposition of the normalized error process $U^N$ as
\begin{equation}
\label{method_AED}
U^{N}_t = Q^N_t + J^N_t + \left( \int_0^t \displaystyle  H^{0,N}_s U^N_s ds + \sum \limits_{j=1}^d \int_0^t \displaystyle H^{j,N}_s U^N_s dW_s^j \right),
\end{equation} 
where $H^{j,N}, j \in  \left\{0,\ldots,d\right\}$ take values in $\mathbb{R}^n \otimes \mathbb{R}^n$, $J^N = -N \bar{J}^N$ is a source term and $Q^N = -N \bar{Q}^N$ is a remainder term. By analyzing the stable convergence in law of the source term $J^N$, we will be able to prove the above result using Theorem 2.5 in \cite{NV2}. 
\subsection{It\^o decomposition of $X^{NV}$}
To get the It\^o decomposition of $X^{NV}$, the main difficulty is to explicit the derivatives of $h_{d+1}$ given by \eqref{eq_hd+1}. Even if the commutativity of the Brownian vector fields $\sigma^j, j \in  \left\{1,\ldots,d\right\}$, simplifies the calculation a lot, it is still cumbersome. To compute the derivatives of $h_{d+1}$, we begin by writing $h_{d+1}$ in integral form as follows. The Frobenius theorem (see \cite{Dieudonee} or \cite{Doss}) ensures that for all permutation $\xi$  on $\left\{1,\ldots,d\right\}$:
\begin{equation}
\begin{split}
\label{h_d+1_Integral}
h_{d+1}\left(t_0,\ldots,t_{d+1};y\right) & = y + \int_{0}^{t_0} \displaystyle \sigma^0 \left( \exp\left(s \sigma^0\right)y\right) ds \\
& + \sum \limits_{k=1}^d \int_{0}^{t_{\xi \left(k\right)}} \displaystyle \sigma^{\xi \left(k\right)} \left( \exp\left(s\sigma^{\xi \left(k\right)}\right) \left( \prod \limits_{m < k} \exp\left(t_{\xi \left(m\right)}\sigma^{\xi \left(m\right)}\right)\right) \exp\left(t_0\sigma^0\right)y \right) ds\\
& +  \int_{0}^{t_{d+1}} \displaystyle \sigma^0 \left( h_{d+1}\left(t_0,\ldots,t_{d},s;y\right) \right) ds.
\end{split}
\end{equation}
\\\\

\textbf{First order derivative with respect to $t_{0}$.}\\\\
To compute the derivative with respect to $t_0$, we choose $\xi = Id$, and we introduce the functions
\begin{equation}
h_j\left(t_0,\ldots,t_j;y\right) = \exp(t_j \sigma^j) \ldots \exp(t_1 \sigma^1) \exp(t_0 \sigma^0)y, 
\end{equation}
for $j \in \left\{0,\ldots,d\right\}$. Then, with $\xi = Id$, \eqref{h_d+1_Integral} becomes
\begin{equation*}
\begin{split}
h_{d+1}\left(t_0,\ldots,t_{d+1};y\right) & = y + \int_{0}^{t_0} \displaystyle \sigma^0 \left( h_0\left(s;y\right)\right) ds  + \sum \limits_{k=1}^d \int_{0}^{t_{k}} \displaystyle \sigma^{k} \left( h_k\left(t_0,\ldots,t_{k-1},s;y\right)\right) ds\\
& +  \int_{0}^{t_{d+1}} \displaystyle \sigma^0 \left( h_{d+1}\left(t_0,\ldots,t_{d},s;y\right) \right) ds,
\end{split}
\end{equation*}
and it follows that
\begin{equation*}
\begin{split}
\partial_{t_0} h_{d+1}\left(t_0,\ldots,t_{d+1};y\right) & = \sigma^0\left(h_0\left(t_0;y\right)\right) +  \sum \limits_{j=1}^d \int_{0}^{t_j} \displaystyle \partial \sigma^j \left( h_j\left(t_0,\ldots,t_{j-1},s;y\right) \right) \partial_{t_0} h_j\left(t_0,\ldots,t_{j-1},s;y\right) ds\\
& +  \int_{0}^{t_{d+1}} \displaystyle \partial \sigma^0 \left( h_{d+1}\left(t_0,\ldots,t_d,s;y\right) \right) \partial_{t_0}h_{d+1}\left(t_0,\ldots,t_d,s;y\right) ds. 
\end{split}
\end{equation*}
Moreover,
\begin{equation*}
\begin{split}
\partial_{t_0} h_{j}\left(t_0,\ldots,t_{j};y\right) & = \sigma^0\left(h_0\left(t_0;y\right)\right) +  \sum \limits_{k=1}^j \int_{0}^{t_k} \displaystyle \partial \sigma^k \left( h_k\left(t_0,\ldots,t_{k-1},s;y\right) \right) \partial_{t_0} h_k\left(t_0,\ldots,t_{k-1},s;y\right) ds.
\end{split}
\end{equation*}
Then, solving these linear differential equations, we obtain  by induction
\begin{equation}
\begin{split}
\partial_{t_0}  h_{d+1}\left(t_0,\ldots,t_{d+1};y\right) & = R^{d+1}\left(t_0,\ldots,t_{d+1};y\right)  R^d\left(t_0,\ldots,t_{d};y\right) \ldots R^1\left(t_0,t_1;y\right) \sigma^0\left(h_0\left(t_0;y\right) \right),
\end{split}
\end{equation}
where, for $j \in \left\{1,\ldots,d+1\right\}, R^j$ is the solution of the following linear system
\begin{equation}
\left\{
    \begin{array}{ll}
\frac{d}{dt}  R^j\left(t_0,\ldots,t_{j-1},t;y\right) = \partial \sigma^j \left(  h_j\left(t_0,\ldots,t_{j-1},t;y\right) \right)   R^j\left(t_0,\ldots,t_{j-1},t;y\right)\\
 R^j\left(t_0,\ldots,t_{j-1},0;y\right) = Id_n,
\end{array}
\right.
\end{equation}
with $\sigma^{d+1} = \sigma^{0}$ by convention.\\\\  
\textbf{First order derivative with respect to $t_{j}, j \in \left\{1,\ldots,d\right\}$.}\\\\
Choosing $\xi$ such that $\xi\left(d\right) = j$, we get
\begin{equation}
\begin{split}
\label{dtjh}
\partial_{t_j} h_{d+1}\left(t_0,\ldots,t_{d+1};y\right) & =  \sigma^j\left(h_d\left(t_0,\ldots,t_d;y\right)\right) +  \int_{0}^{t_{d+1}} \displaystyle \partial \sigma^0 \left(  h_{d+1}\left(t_0,\ldots,t_d,s;y\right) \right) \partial_{t_j}  h_{d+1}\left(t_0,\ldots,t_d,s;y\right)  ds.
\end{split}
\end{equation}
Solving this linear differential equation, we obtain
\begin{equation}
\partial_{t_j} h_{d+1}\left(t_0,\ldots,t_{d+1};y\right) = R^{d+1}\left(t_0,\ldots,t_{d+1};y\right) \sigma^j\left(h_d\left(t_0,t_1,\ldots,t_d;y\right)\right).
\end{equation}

\textbf{First order derivative with respect to $t_{d+1}$.}\\\\
The derivative with respect to $t_{d+1}$ is trivial:
\begin{equation}
\partial_{t_{d+1}} h_{d+1}\left(t_0,\ldots,t_{d+1};y\right) = \sigma^0\left(h_{d+1}\left(t_0,\ldots,t_{d+1};y\right)\right).
\end{equation}

\textbf{Second order derivative with respect to $t_{j}, j \in \left\{1,\ldots,d\right\}$.}\\\\
Using \eqref{dtjh} and  
\begin{equation}
\partial_{t_j} h_{d}\left(t_0,\ldots,t_{d};y\right)  =  \sigma^j\left(h_d\left(t_0,\ldots,t_d;y\right)\right), 
\end{equation}
we have
\begin{equation*}
\begin{split}
\partial^2_{t_jt_j} h_{d+1}\left(t_0,\ldots,t_{d+1};y\right) &  = \partial \sigma^j \sigma^j \left(h_d\left(t_0,\ldots,t_d;y\right)\right) \\
& + \int_{0}^{t_{d+1}} \displaystyle \left( \partial^2 \sigma^0 \left(  h_{d+1}\left(t_0,\ldots,t_d,s;y\right) \right)\odot \partial_{t_j} h_{d+1}\left(t_0,\ldots,t_d,s;y\right) \right) \partial_{t_j}  h_{d+1}\left(t_0,\ldots,t_d,s;y\right)  ds  \\
& + \int_{0}^{t_{d+1}} \displaystyle \partial \sigma^0 \left(  h_{d+1}\left(t_0,\ldots,t_d,s;y\right) \right) \partial^2_{t_jt_j}  h_{d+1}\left(t_0,\ldots,t_d,s;y\right) ds.
\end{split}
\end{equation*}
Solving this linear differential equation, we get
\begin{multline}
\partial^2_{t_jt_j} h_{d+1}\left(t_0,\ldots,t_{d+1};y\right) = R^{d+1}\left(t_0,\ldots,t_{d+1};y\right)  \partial \sigma^j \sigma^j \left(h_d\left(t_0,\ldots,t_d;y\right)\right) \\
+   \int_{0}^{t_{d+1}} \displaystyle R^{d+1}\left(t_0,\ldots,t_{d+1};y\right)\left(R^{d+1}\right)^{-1}\left(t_0,\ldots,t_d,s;y\right)\left( \left( \partial^2 \sigma^0 \circ  h_{d+1} \odot \partial_{t_j} h_{d+1} \right) \partial_{t_j}  h_{d+1}\right) \left(t_0,\ldots,t_d,s;y\right)  ds.
\end{multline}
\textbf{It\^o's formula.}\\\\
Using It\^o's formula, we obtain:
\begin{equation}
\begin{split}
\label{dXNV0}
d X^{NV}_t &= \frac{1}{2} \partial_{t_{0}}   h_{d+1}\left(\frac{\Delta t}{2}, \Delta W_t, \frac{\Delta t}{2}; X^{NV}_{\hat{\tau}_{t}} \right) dt + \sum \limits_{j=1}^d  \partial_{t_j}  h_{d+1}\left(\frac{\Delta t}{2}, \Delta W_t, \frac{\Delta t}{2}; X^{NV}_{\hat{\tau}_{t}} \right) dW^j_t\\
&+  \frac{1}{2} \partial_{t_{d+1}}  h_{d+1}\left(\frac{\Delta t}{2} , \Delta W_t, \frac{\Delta t}{2}; X^{NV}_{\hat{\tau}_{t}} \right)  dt  + \frac12 \sum \limits_{j=1}^d  \partial^2_{t_j t_j}  h_{d+1}\left(\frac{\Delta t}{2}, \Delta W_t, \frac{\Delta t}{2}; X^{NV}_{\hat{\tau}_{t}} \right) dt\\
&= \frac{1}{2} Z^0_t dt + \sum \limits_{j=1}^d   Y^{d+1}_{t,t} \sigma^j\left(\bar{X}_{t,t}\right) dW^j_t+ \frac{1}{2}  \sigma^0\left(X^{NV}_t\right)  dt + \frac12 \sum \limits_{j=1}^d    \left( Y^{d+1}_{t,t} \partial \sigma^j \sigma^j\left(\bar{X}_{t,t}\right) +  Z^{d+1,j}_t \right)dt, 
\end{split} 
\end{equation}
where for $t \in [0,T]$
\begin{equation*}
Z^0_t =  Y^{d+1}_{t,t} \ldots  Y^{1}_{t,t} \sigma^0\left(\bar{X}_t^0 \right), 
\end{equation*}
for $j \in \left\{1,\ldots,d\right\}, s \in [\hat{\tau}_{t},t]$
\begin{equation*}
 Y^{j}_{t,s} =  R^{j}\left(\frac{\Delta t}{2}, \Delta W^1_t,\ldots,\Delta W^{j-1}_t,\Delta W^{j}_s  ; X^{NV}_{\hat{\tau}_{t}} \right), 
\end{equation*}
\begin{equation*}
 Y^{d+1}_{t,s} =  R^{d+1}\left(\frac{\Delta t}{2}, \Delta W^1_t,\ldots,\Delta W^{d}_t,\frac{\Delta s}{2}  ; X^{NV}_{\hat{\tau}_{t}} \right), 
\end{equation*}
\begin{equation*}
\bar{X}^0_{t} =  h_0\left(\frac{\Delta t}{2};X^{NV}_{\hat{\tau}_{t}}\right),
\end{equation*}
\begin{equation*}
\bar{X}_{t,s} =  h_d\left(\frac{\Delta t}{2},\Delta W_s;X^{NV}_{\hat{\tau}_{t}}\right),
\end{equation*}
\begin{equation*}
Z^{d+1,j}_t = \frac{1}{2}\int_{\hat{\tau}_{t}}^{t} \displaystyle  Y_{t,t}^{d+1} \left(Y^{d+1}_{t,s}\right)^{-1} \displaystyle \left( \partial^2 \sigma^0\left(\bar{X}^{d+1}_{t,s}\right) \odot \left(Y^{d+1}_{t,s} \sigma^j\left(\bar{X}_{t,t}\right)\right)\right) Y^{d+1}_{t,s} \sigma^j\left(\bar{X}_{t,t}\right)\ ds,
\end{equation*}
and
\begin{equation}
\bar{X}^{d+1}_{t,s} =  h_{d+1}\left(\frac{\Delta t}{2},\Delta W_t,\frac{\Delta s}{2};X^{NV}_{\hat{\tau}_{t}}\right).
\end{equation}
The next lemma gives estimations of the moment of the Ninomiya-Victoir scheme $X^{NV}$ and its increments. It also compares $X^{NV}$ to the intermediate processes $\bar{X}^0$, $\bar{X}$ and $\bar{X}^{d+1}$.
This result is very similar to Lemmas 4.4 and 4.5 in \cite{NV2}, that is why we omit its proof.
\begin{alem}
\label{Lemme0}
Assume that 
\begin{itemize}
\item $\forall j \in \left\{1,\ldots,d\right\}, \sigma^j \in \mathcal{C}^{1}\left(\mathbb{R}^n,\mathbb{R}^n\right)$  with bounded first order derivatives,
\item $\sigma^0$ and $\sum \limits_{j=1}^d \partial \sigma^j \sigma^j$ are Lipschitz continuous functions.
\end{itemize}
Then, $\forall p \ge 1, \exists C_0 \in \mathbb{R}_+^*,  \forall N \in \mathbb{N}^*, \forall t \in [0,T],$
\begin{equation}
\begin{split}
\label{NV_moment}
\mathbb{E}\left[ \left\|  X^{NV}_{t} \right\|^{2p}\right] & \leq   C_0,
\end{split}
\end{equation}
\begin{equation}
\begin{split}
\label{Acc_NV}
\mathbb{E}\left[ \left\| X^{NV}_{t}  - X^{NV}_{\hat{\tau}_t}\right\|^{2p}\right] & \leq   C_0 \left(\Delta t\right)^{p},
\end{split}
\end{equation}
\begin{equation}
\begin{split}
\label{Acc_0}
\mathbb{E}\left[ \left\| X^{NV}_{t}  - \bar{X}^0_{s}\right\|^{2p}\right] & \leq   C_0 \left(\Delta t\right)^{p}, \forall s \in [\hat{\tau}_t,t],
\end{split}
\end{equation}
\begin{equation}
\begin{split}
\label{Acc_}
\mathbb{E}\left[ \left\| X^{NV}_{t}  - \bar{X}_{t,s}\right\|^{2p}\right] & \leq   C_0 \left(\Delta t\right)^{p},\forall s \in [\hat{\tau}_t,t],
\end{split}
\end{equation}
\begin{equation}
\begin{split}
\label{Acc_d+1}
\mathbb{E}\left[ \left\| X^{NV}_{t}  - \bar{X}^{d+1}_{t,s}\right\|^{2p}\right] & \leq   C_0 \left(\Delta t\right)^{2p}, \forall s \in [\hat{\tau}_t,t].
\end{split}
\end{equation}
Moreover if the commutativity condition \eqref{Comm_condition} holds:
\begin{equation}
\begin{split}
\label{Acc_t}
\mathbb{E}\left[ \left\| X^{NV}_{t}  - \bar{X}_{t,t}\right\|^{2p}\right] & \leq   C_0 \left(\Delta t\right)^{2p}.
\end{split}
\end{equation}
\end{alem}

\subsection{Suitable decomposition of $X^{NV}$}

Our derivation of a decomposition of the form \eqref{decomposition_1_XVN} is carried out in two steps. The first step will consist in approximating the dynamics of $X^{NV}$ with strong order $3/2$. In the second step we will identify the appropriate source term $J^N$.
\subsubsection{Approximation with strong order $3/2$}

The dynamics \eqref{dXNV0} of $X^{NV}$ is not really tractable to study the normalized error process $U^N$.
In the following, we provide an approximation with strong order $3/2$ of the theoretical dynamics of $X^{NV}$. 
The goal is to be able to write the dynamics of $X^{NV}$ in the form \eqref{decomposition_1_XVN}.
We begin by indicating a natural approximation with strong order $2$, respectively $3/2$, of the intermediate processes $Y^{d+1}$ and $\left(\sigma^j\left(\bar{X}_{t,t}\right)\right)_{0\leq t \leq T}$, respectively $\left(\partial \sigma^j  \sigma^j\left(\bar{X}_{t,t}\right)\right)_{0\leq t \leq T}$ and $Z^{d+1,j}$. \\\\
\\\\\

\begin{aprop}
\label{Prop_App_1}
Assume that 
\begin{itemize}
\item $\forall j \in \left\{0,\ldots,d\right\}, \sigma^j \in \mathcal{C}^{2}\left(\mathbb{R}^n,\mathbb{R}^n\right)$  with bounded first order derivatives, $\partial \sigma^j \in LIP_{loc}^{pgc}\left(\mathbb{R}^n \otimes \mathbb{R}^n\right)$ and $\partial^2 \sigma^j  \in LIP_{loc}^{pgc}\left(\mathbb{R}^n \otimes \mathbb{R}^n \otimes \mathbb{R}^n\right)$,
\item $\sum \limits_{j=1}^d \partial \sigma^j \sigma^j$ is a Lipschitz continuous function.
\end{itemize}
Then, $\forall p \ge 1, \exists C_1 \in \mathbb{R}_+^*, \forall N \in \mathbb{N}^*, \forall t \in [0,T], \forall j \in \left\{1,\ldots,d\right\},$
\begin{equation}
\begin{split}
\label{approx_1}
\mathbb{E}\left[ \left\| Y^{d+1}_{t,s}  - Id_n - \frac{1}{2} \Delta s  \partial \sigma^0\left(X^{NV}_t\right) \right\|^{2p}\right] & \leq   C_1 \left(\Delta t\right)^{4p}, \forall s \in [\hat{\tau}_{t},t],
\end{split}
\end{equation}
\begin{equation}
\begin{split}
\label{approx_4}
\mathbb{E}\left[ \left\|  Z^{d+1,j}_{t} -  \frac{1}{2} \Delta t \left( \partial^2\sigma^0\odot \sigma^j\right) \sigma^j\left(X^{NV}_{\hat{\tau}_t}\right)\right\|^{2p}\right] \leq   C_1 \left(\Delta t\right)^{3p}.
\end{split}
\end{equation}
Moreover if the commutativity condition \eqref{Comm_condition} holds:
\begin{equation}
\label{approx_2}
\begin{split}
\mathbb{E}\left[ \left\| \sigma^j\left(\bar{X}_{t,t}\right) - \sigma^j\left(X^{NV}_{t}\right) + \frac{1}{2} \Delta t \partial \sigma^j \sigma^0 \left(X^{NV}_{t} \right)\right\|^{2p}\right] \leq   C_1 \left(\Delta t\right)^{4p},
\end{split}
\end{equation}
\begin{equation}
\begin{split}
\label{approx_3}
\mathbb{E}\left[ \left\| \partial\sigma^j\sigma^j\left(\bar{X}_{t,t}\right) - \partial\sigma^j\sigma^j\left(X^{NV}_{t}\right) + \frac{1}{2} \Delta t \left(\partial^2\sigma^j \odot \sigma^j  + \left(\partial \sigma^j\right)^2 \right) \sigma^0 \left(X^{NV}_{\hat{\tau}_t} \right)\right\|^{2p}\right] \leq   C_1 \left(\Delta t\right)^{3p}.
\end{split}
\end{equation}
\end{aprop}
In \eqref{approx_4} and \eqref{approx_3}, we have replaced by $X^{NV}_{\hat{\tau}_t}$ the argument of the functions which would naturally appears in Taylor expansion when Lemma \ref{Lemme0} ensures that the strong order $3/2$ is preserved.
Although the above approximations are very intuitive, their proofs are both heavy and technical.
That is why, the proof of this proposition is postponed to the Appendix.
To obtain an approximation with strong order $3/2$ of the form \eqref{decomposition_1_XVN}, it remains to estimate $\displaystyle \int_{0}^{t} Z^0_s ds$, for $t \in [0,T]$.
\begin{aprop}
\label{Prop_App_2}
Assume that 
\begin{itemize}
\item $ \sigma^0 \in \mathcal{C}^{1}\left(\mathbb{R}^n,\mathbb{R}^n\right)$  with bounded first order derivatives, and $\partial \sigma^0 \in LIP_{loc}^{pgc}\left(\mathbb{R}^n \otimes \mathbb{R}^n\right)$,
\item $\forall j \in \left\{1,\ldots,d\right\},\sigma^j \in \mathcal{C}^{2}\left(\mathbb{R}^n,\mathbb{R}^n\right)$  with bounded first order derivatives, $\partial \sigma^j \in LIP_{loc}^{pgc}\left(\mathbb{R}^n \otimes \mathbb{R}^n\right)$, and $\partial^2 \sigma^j \in LIP_{loc}^{pgc}\left(\mathbb{R}^n \otimes \mathbb{R}^n \otimes \mathbb{R}^n\right)$,
\item $\forall j \in \left\{1,\ldots,d\right\}, \partial \sigma^j \sigma^j$ is a Lipschitz continuous function.
\end{itemize}
 Then, denoting by, 
\begin{equation}
\begin{split}
\label{sub_theta}
\theta_t^0 &=   \sigma^0\left(\bar{X}^0_t\right) + \sum \limits_{j=1}^d  \Delta W^j_t \partial \sigma^j \sigma^0\left(X^{NV}_{\hat{\tau}_t}\right)  + \frac{1}{2} \sum \limits_{j=1}^d \Delta t\left(\partial^2 \sigma^j \odot \sigma^j + \left(\partial \sigma^j\right)^2\right) \sigma^0\left(X^{NV}_{\hat{\tau}_t}\right)  + \frac{1 }{2} \Delta t \partial \sigma^0 \sigma^0\left(X^{NV}_{\hat{\tau}_t}\right),
\end{split}
\end{equation}
for $t \in [0,T]$ ,we have: 
\begin{equation}
\begin{split}
\label{approx_5}
\forall p \ge 1, \exists C_2 \in \mathbb{R}_+^*, \forall N \in \mathbb{N}^*, \mathbb{E}\left[ \underset{t\leq T}{\sup} \left\| \displaystyle \int_{0}^{t} \left( Z^0_s - \theta_s^0\right) ds \right\|^{2p}\right] \leq C_2 h^{3p}.
\end{split}
\end{equation}
\end{aprop}
This approximation is not really intuitive. Indeed, to get this result, we approximate the process $Z^0$, and then we use the integration by parts formula to identify the dominant contribution.  The proof of this proposition is also postponed to the Appendix.

We are now able to approximate the dynamics of the Ninomiya-Victoir scheme. 
Using $\displaystyle b = \sigma^0 + \frac12 \sum \limits_{j=1}^d \partial \sigma^j  \sigma^j$, together with Propositions \ref{Prop_App_1} and \ref{Prop_App_2}, it is easy to see that an adequate decomposition of the dynamics of $X^{NV}$ is given by
\begin{equation}
\begin{split}
\label{CC_Approx}
d X^{NV}_t &=  b\left(X^{NV}_t\right) dt + \sum \limits_{j=1}^d  \sigma^j\left(X^{NV}_t\right) dW^j_t +  \frac12  \left(\sigma^0\left(\bar{X}^0_t\right) - \sigma^0\left(X^{NV}_t\right) \right) dt\\
& + \frac{1}{2} \sum \limits_{j=1}^d  \Delta t \left[ \sigma^j, \sigma^0\right ]\left(X^{NV}_t\right) dW^j_t + \frac12 \sum \limits_{j=1}^d  \Delta W^j_t \partial \sigma^j \sigma^0\left(X^{NV}_{\hat{\tau}_t}\right) dt \\
& +\frac{1}{4} \Delta t \left( \partial \sigma^0 \sigma^0 + \sum \limits_{j=1}^d \left( \partial^2\sigma^0 \odot \sigma^j + \partial\sigma^0\partial\sigma^j  \right)\sigma^j  \right)\left(X^{NV}_{\hat{\tau}_t}\right) dt + d\bar{Q}^{1,N}_t,
\end{split} 
\end{equation}
where, $\bar{Q}^{1,N}$ is defined by $\bar{Q}^{1,N}_0 = 0$ and 
\begin{equation}
\begin{split}
\label{R1}
d \bar{Q}^{1,N}_t &=  \frac12 \left( Z^0_t - \theta^0_t \right) + \sum \limits_{j=1}^d \left( Y^{d+1}_{t,t} \sigma^j\left(\bar{X}_{t,t}\right) - \sigma^j\left(X^{NV}_t\right) - \frac12 \Delta t \left[ \sigma^j, \sigma^0\right ]\left(X^{NV}_t\right) \right) dW^j_t\\
& + \frac12 \sum \limits_{j=1}^d    \left( Y^{d+1}_{t,t} \partial \sigma^j \sigma^j\left(\bar{X}_{t,t}\right) - \partial\sigma^j \sigma^j\left(X^{NV}_t\right)  - \frac12 \Delta t \left( \partial \sigma^0 \partial \sigma^j \sigma^j - \left(\partial^2 \sigma^j \odot \sigma^j + \left(\partial \sigma^j\right)^2\right) \sigma^0 \right)\left(X^{NV}_{\hat{\tau}_t}\right) \right)  dt  \\
& + \frac12 \sum \limits_{j=1}^d  \left( Z^{d+1,j}_t - \frac12 \Delta t \left( \partial^2\sigma^0 \odot \sigma^j \right)\sigma^j\left(X^{NV}_{\hat{\tau}_t}\right) \right) dt,
\end{split} 
\end{equation} 
is a remainder term with strong order $3/2$. The proof of the following proposition is also postponed to the Appendix.
\begin{aprop}
\label{Prop_Q}
Under the assumptions of Theorem \ref{T_EP_CC}:
 $\forall p \ge 1, \exists C_3 \in \mathbb{R}_+^*, \forall N \in \mathbb{N}^*,$ 
\begin{equation}
\begin{split}
\label{Prop_est}
\mathbb{E}\left[ \underset{t\leq T}{\sup} \left\| \bar{Q}^{1,N}_t  \right\|^{2p}\right] & \leq   C_3 h^{3p}.
\end{split}
\end{equation}
\end{aprop}

\subsubsection{Identification of the source term $J^N$}
The goal of this subsection is to identify the source term $J^N$. Deducing the strong convergence with order $1$ of the Ninomiya-Victoir scheme from the equation \eqref{CC_Approx} does not look straightforward.
Indeed, at first glance, in \eqref{CC_Approx}, the terms 
\begin{equation}
\label{IPP1}
\frac{N}2 \displaystyle \int_{0}^T   \left( \sigma^0\left(X^{NV}_t\right) - \sigma^0\left(\bar{X}^0_t\right) \right)dt
\end{equation}
and 
\begin{equation}
\label{IPP2}
\frac{N}2 \sum \limits_{j=1}^d  \displaystyle \int_{0}^T     \Delta W^j_t \partial \sigma^j \sigma^0\left(X^{NV}_{\hat{\tau}_t}\right)dt
\end{equation}
seem to diverge as $N$ goes to infinity. Actually, using the integration by parts formula, we show that both terms, \eqref{IPP1} and \eqref{IPP2}, are bounded by a constant independent of $N$ in $L^2$, which is consistent with Theorem \ref{SC2}. More precisely, we have the following result. 
\begin{aprop}
\label{LemmeIPP}
Let 
\begin{equation}
\begin{split}
I^0_t  &= \sum \limits_{j=1}^d  \int_{0}^{t} \displaystyle \left(t \wedge  \check{\tau}_s - s\right) \partial \sigma^0 \sigma^j \left(X^{NV}_s\right) dW^j_s  
-    \displaystyle \frac 12 \int_{0}^{t} \left(t \wedge  \check{\tau}_s - s\right) \partial \sigma^0 \sigma^0\left( X^{NV}_{\hat{\tau}_s} \right) ds\\
& - \displaystyle  \frac 12  \sum \limits_{j=1}^d  \int_{0}^{t} \left(t \wedge \check{\tau}_s - s\right) \left(\partial^2 \sigma^0 \odot \sigma^j + \partial \sigma^0\partial \sigma^j  \right) \sigma^j\left(X^{NV}_{\hat{\tau}_s}\right) ds,~ t \in [0,T],
\end{split}
\end{equation}
 and for $j \in \left\{1,\ldots,d\right\}$
\begin{equation}
I^j_t  = \int_{0}^{t} \displaystyle \left(t \wedge \check{\tau}_s - s\right)\partial \sigma^j\sigma^0\left(X^{NV}_{s}  \right) dW^j_s,~t \in [0,T].
\end{equation}
Under the assumptions of Theorem \ref{T_EP_CC}: $\forall p \ge 1, \exists C_4 \in \mathbb{R}_+^*, \forall N \in \mathbb{N}^*, \forall j \in \left\{1,\ldots,d\right\}$:
\begin{equation}
\begin{split}
\label{IPP1_2}
\mathbb{E}\Bigg[ \underset{t\leq T}{\sup} \Bigg\| & \displaystyle  \int_{0}^{t} \displaystyle \sigma^0\left(X^{NV}_s\right)  - \sigma^0\left(\bar{X}^0_{s}  \right) ds  -  I^0_t \Bigg\|^{2p}\Bigg]  \leq   C_4 h^{3p},
\end{split} 
\end{equation} 
\begin{equation}
\begin{split}
\label{IPP1_1}
\mathbb{E}\left[ \underset{t\leq T}{\sup} \left\| \displaystyle \int_{0}^t   \Delta W^j_s \partial \sigma^j \sigma^0\left(X^{NV}_{\hat{\tau}_s}\right)ds - I_t^j \right\|^{2p}\right] & \leq   C_4 h^{3p}.
\end{split}
\end{equation}
\end{aprop}
\begin{adem}
We start by proving \eqref{IPP1_1} and we denote by 
\begin{equation}
\Phi^{N}_t = \int_{0}^{t} \displaystyle  \Delta W^j_s  \partial\sigma^j\sigma^0\left(X^{NV}_{\hat{\tau}_s}  \right) ds - I^j_t.
\end{equation}
Using the integration by parts formula, we have
\begin{equation*}
\int_{0}^{t} \displaystyle \Delta W^j_s \partial\sigma^j\sigma^0\left(X^{NV}_{\hat{\tau}_s}  \right)  ds  = \int_{0}^{t} \displaystyle  \left(t \wedge  \check{\tau}_s - s\right) \partial\sigma^j\sigma^0\left(X^{NV}_{\hat{\tau}_s}  \right) d W^j_s.   
\end{equation*}
Therefore, 
\begin{equation*}
\begin{split}
\Phi^{N}_t &=  - \int_{0}^{t} \displaystyle \left(t \wedge \check{\tau}_s - s\right)\left(\partial \sigma^j\sigma^0\left(X^{NV}_{s}\right) - \partial \sigma^j\sigma^0\left(X^{NV}_{\hat{\tau}_s}  \right) \right) dW^j_s.
\end{split} 
\end{equation*} 
Combining the Burkholder-Davis-Gundy inequality and a convexity inequality, we get a constant $\alpha_1$ independent of $N$ such that
\begin{equation}
\begin{split}
\mathbb{E}\left[ \underset{t\leq T}{\sup} \left\| \Phi^{N}_t  \right\|^{2p}\right] & \leq   \alpha_1 \int_{0}^{T} \displaystyle \left(t \wedge \check{\tau}_s - s\right)^{2p} \mathbb{E}\left[ \left\|  \partial \sigma^j\sigma^0\left(X^{NV}_{s}\right) - \partial \sigma^j\sigma^0\left(X^{NV}_{\hat{\tau}_s}  \right)  \right\|^{2p}\right] ds.
\end{split}
\end{equation}
As the function $\partial \sigma^j\sigma^0 \in LIP_{loc}^{pgc}\left(\mathbb{R}^n \right)$ is locally Lipschitz with polynomially growing Lipschitz constant, there exists $c \in \mathbb{R}_+^*$ and $q \in \mathbb{N}$ such that
\begin{equation}
\begin{split}  
\left\| \partial \sigma^j\sigma^0\left(X^{NV}_{s}\right) - \partial \sigma^j\sigma^0\left(X^{NV}_{\hat{\tau}_s}  \right)  \right\|  \leq c \left(1 + \left\| X^{NV}_s \right\|^q \vee \left\| X^{NV}_{\hat{\tau}_s} \right\|^q\right) \left\| X^{NV}_{s} - X^{NV}_{\hat{\tau}_s} \right\|.    
\end{split}
\end{equation}
Hence, 
\begin{equation}
\begin{split}  
\mathbb{E}\left[ \left\| \partial \sigma^j\sigma^0\left(X^{NV}_{s}\right) - \partial \sigma^j\sigma^0\left(X^{NV}_{\hat{\tau}_s}  \right) \right\|^{2p}\right]  \leq c^{2p}  \mathbb{E}\left[  \left(1 + \left\| X^{NV}_s \right\|^q \vee \left\| X^{NV}_{\hat{\tau}_s} \right\|^q\right)^{2p} \left\| X^{NV}_{s} - X^{NV}_{\hat{\tau}_s}    \right\|^{2p}\right]. 
\end{split}
\end{equation}
Applying the Cauchy-Schwarz inequality together with \eqref{NV_moment} and \eqref{Acc_NV} from Lemma \ref{Lemme0}, we easily get a constant $\alpha_2 \in \mathbb{R}_+^*$ independent of $N$ such that:  
\begin{equation}
\begin{split}
\mathbb{E}\left[ \left\| \partial \sigma^j\sigma^0\left(X^{NV}_{s}\right) - \partial \sigma^j\sigma^0\left(X^{NV}_{\hat{\tau}_s}  \right) \right\|^{2p}\right] & \leq   \alpha_2 h^{p},
\end{split}
\end{equation}
and we conclude that
\begin{equation}
\begin{split}
\mathbb{E}\left[ \underset{t\leq T}{\sup} \left\| \Phi^{N}_t  \right\|^{2p}\right] & \leq   \alpha_1 \alpha_2 T h^{3p}. 
\end{split}
\end{equation}
Now, we focus on \eqref{IPP1_2} and we denote
\begin{equation}
\begin{split}
\Psi^{N}_t &= \displaystyle  \int_{0}^{t} \displaystyle \left(\sigma^0\left(X^{NV}_s\right)  - \sigma^0\left(\bar{X}^0_{s}  \right) \right)ds  -  I^0_t.
\end{split}
\end{equation}
To get a clearer picture, the $i$-th coordinate, $i \in \left\{1,\ldots,n\right\}$, $\Psi^{i,N}_t$ is given by
\begin{equation}
\begin{split}
\Psi^{i,N}_t &= \displaystyle  \int_{0}^{t} \displaystyle \sigma^{i0}\left(X^{NV}_s\right)  - \sigma^{i0}\left(\bar{X}^0_{s}  \right) ds  -  \sum \limits_{j=1}^d  \int_{0}^{t} \displaystyle \left(t \wedge  \check{\tau}_s - s\right) \nabla \sigma^{i0}\left(X^{NV}_s\right) . ~\sigma^j \left(X^{NV}_s\right) dW^j_s    \\
& -   \displaystyle \frac 12 \int_{0}^{t}  \left(t \wedge  \check{\tau}_s - s\right) \nabla \sigma^{i0}\left( X^{NV}_{\hat{\tau}_s} \right) . ~\sigma^0\left( X^{NV}_{\hat{\tau}_s} \right) ds - \displaystyle  \frac 12 \sum \limits_{j=1}^d \int_{0}^{t} \left(t \wedge \check{\tau}_s - s\right)\nabla \sigma^{i0}\left( X^{NV}_{\hat{\tau}_s} \right) .\left(\partial\sigma^j\sigma^j\left( X^{NV}_{\hat{\tau}_s} \right)    \right) ~ ds   \\
& - \displaystyle  \frac 12 \sum \limits_{j=1}^d \sum \limits_{k=1}^n \sum \limits_{l=1}^n \int_{0}^{t} \left(t \wedge \check{\tau}_s - s\right) \partial^2_{x_k x_l} \sigma^{i0}\left(X^{NV}_{\hat{\tau}_s}\right) \sigma^{kj}\left(X^{NV}_{\hat{\tau}_s}\right) \sigma^{lj}\left(X^{NV}_{\hat{\tau}_s}\right) ds.
\end{split}
\end{equation}
Using once again the integration by parts formula, since $\forall s \in [0,T], X^{NV}_{\hat{\tau}_s} =\bar{X}^0_{\hat{\tau}_s}$, we have
\begin{equation*}
\begin{split}
\int_{0}^{t} \displaystyle \left(\sigma^{i0}\left(X^{NV}_s\right)  - \sigma^{i0}\left(\bar{X}^0_{s}  \right) \right)ds &=\int_{0}^{t} \displaystyle \left(t \wedge \check{\tau}_s - s\right) d\left( \sigma^{i0}\left(X^{NV}_s\right)  - \sigma^{i0}\left(\bar{X}^0_{s}  \right) \right).
\end{split} 
\end{equation*}
Before applying It\^o's formula, we recall that the dynamics of $\bar{X}^{0}$ is given by
\begin{equation}
\begin{split}
d \bar{X}^{0}_t &=  \frac12  \sigma^0\left(\bar{X}^0_t\right) dt,
\end{split} 
\end{equation}  
and that the dynamics of $X^{NV}$  is given by \eqref{CC_Approx}. Since $ \forall t \in [0,T], \forall s \leq t, \left |t \wedge \check{\tau}_s - s\right|\ \leq h$, we rewrite the dynamics of $X^{NV}$ as the sum of its dominant contribution and a remainder term with strong order $1/2$. Using \eqref{CC_Approx} together with $b = \sigma^0 + \displaystyle \frac{1}{2} \sum \limits_{j=1}^d \partial \sigma^j \sigma^j $, we obtain
\begin{equation}
\begin{split}
d X^{NV}_t &=  \frac{1}{2} \sigma^0\left(X^{NV}_t\right) dt + \sum \limits_{j=1}^d  \sigma^j\left(X^{NV}_t\right) dW^j_t + \frac{1}{2} \sum \limits_{j=1}^d \partial\sigma^j \sigma^j\left(X^{NV}_t\right) dt +  \frac12  \sigma^0\left(\bar{X}^0_t\right) dt + d\vartheta_t 
\end{split} 
\end{equation}  
where 
\begin{equation}
\begin{split}
 d\vartheta_t  &=  \frac{1}{2} \sum \limits_{j=1}^d  \Delta t \left[ \sigma^j, \sigma^0\right ]\left(X^{NV}_t\right) dW^j_t + \frac12 \sum \limits_{j=1}^d  \Delta W^j_t \partial \sigma^j \sigma^0\left(X^{NV}_{\hat{\tau}_t}\right) dt \\
& +\frac{1}{4} \Delta t  \left( \partial \sigma^0 \sigma^0 + \sum \limits_{j=1}^d \left( \partial^2\sigma^0 \odot \sigma^j + \partial\sigma^0\partial\sigma^j  \right)\sigma^j  \right)\left(X^{NV}_{\hat{\tau}_t}\right) dt + d\bar{Q}^{1,N}_t.
\end{split} 
\end{equation}  
Notice that $\vartheta$  is a term with strong order $1$ since we have already proved \eqref{IPP1_1}.
Applying It\^o's formula we get
\begin{equation*}
\begin{split}
\int_{0}^{t} \displaystyle \sigma^{i0}\left(X^{NV}_s\right)  - \sigma^{i0}\left(\bar{X}^0_{s}  \right) ds &= \displaystyle \frac12 \int_{0}^{t}  \left(t \wedge \check{\tau}_s - s\right) \nabla \sigma^{i0}\left(X^{NV}_s\right) .~ \sigma^0\left(X^{NV}_s\right) ds \\
& +\sum \limits_{j=1}^d  \int_{0}^{t} \displaystyle \left(t \wedge  \check{\tau}_s - s\right) \nabla \sigma^{i0}\left(X^{NV}_s\right) . ~\sigma^j \left(X^{NV}_s\right) dW^j_s\\
& + \displaystyle \frac12 \sum \limits_{j=1}^d  \int_{0}^{t} \displaystyle \left(t \wedge  \check{\tau}_s - s\right) \nabla \sigma^{i0}\left(X^{NV}_s\right) . ~\partial\sigma^j  \sigma^j \left(X^{NV}_s\right) ds\\
& +\displaystyle \frac12 \int_{0}^{t}  \left(t \wedge \check{\tau}_s - s\right)  \left(\nabla \sigma^{i0}\left(X^{NV}_s\right) - \nabla \sigma^{i0}\left(\bar{X}^0_{s}\right)\right) .~ \sigma^0\left(\bar{X}^0_{s}\right)  ds \\
& + \frac12 \sum \limits_{j=1}^d \sum \limits_{k=1}^n \sum \limits_{l=1}^n \int_{0}^{t} \displaystyle \left(t \wedge \check{\tau}_s - s\right) \partial^2_{x_k x_l} \sigma^{i0}\left(X^{NV}_s\right) \sigma^{kj}\left(X^{NV}_s\right) \sigma^{lj}\left(X^{NV}_s\right) ds\\
& +\displaystyle \int_{0}^{t}  \left(t \wedge \check{\tau}_s - s\right)  \nabla \sigma^{i0}\left(X^{NV}_s\right).~d\vartheta_s \\
& + \frac12 \sum \limits_{k=1}^n \sum \limits_{l=1}^n \int_{0}^{t} \displaystyle \left(t \wedge \check{\tau}_s - s\right) \partial^2_{x_k x_l} \sigma^{i0}\left(X^{NV}_s\right) d\langle \vartheta^k,\vartheta^l \rangle_s\\
&+  \frac12 \sum \limits_{j=1}^d \sum \limits_{k=1}^n \sum \limits_{l=1}^n \int_{0}^{t} \displaystyle \left(t \wedge \check{\tau}_s - s\right) \partial^2_{x_k x_l} \sigma^{i0}\left(X^{NV}_s\right) \sigma^{lj}\left(X^{NV}_s\right) d\langle \vartheta^k,W^j \rangle_s.
\end{split}
\end{equation*}
Then, it follows that
\begin{equation*}
\begin{split}
\Psi^{i,N}_t &=  \displaystyle \frac12 \int_{0}^{t}  \left(t \wedge \check{\tau}_s - s\right) \left(\nabla \sigma^{i0}\left(X^{NV}_s\right) .~ \sigma^0\left(X^{NV}_s\right) - \nabla \sigma^{i0}\left(X^{NV}_{\hat{\tau}_s}\right) .~ \sigma^0\left(X^{NV}_{\hat{\tau}_s}\right)\right)ds \\
& + \displaystyle \frac12 \sum \limits_{j=1}^d  \int_{0}^{t} \displaystyle \left(t \wedge  \check{\tau}_s - s\right)\left( \nabla \sigma^{i0}\left(X^{NV}_s\right) . ~\partial\sigma^j  \sigma^j \left(X^{NV}_s\right) -  \nabla \sigma^{i0}\left(X^{NV}_{\hat{\tau}_s}\right) . ~\partial\sigma^j  \sigma^j \left(X^{NV}_{\hat{\tau}_s}\right)\right) ds\\
& +\displaystyle \frac12 \int_{0}^{t}  \left(t \wedge \check{\tau}_s - s\right)  \left(\nabla \sigma^{i0}\left(X^{NV}_s\right) - \nabla \sigma^{i0}\left(\bar{X}^0_{s}\right)\right) .~ \sigma^0\left(\bar{X}^0_{s}\right)  ds \\
& + \frac12 \sum \limits_{j=1}^d \sum \limits_{k=1}^n \sum \limits_{l=1}^n \int_{0}^{t} \displaystyle \left(t \wedge \check{\tau}_s - s\right) \left(\partial^2_{x_k x_l} \sigma^{i0} \sigma^{kj} \sigma^{lj}\left(X^{NV}_s\right) - \partial^2_{x_k x_l} \sigma^{i0} \sigma^{kj} \sigma^{lj}\left(X^{NV}_{\hat{\tau}_s}\right)\right) ds\\
& +\displaystyle \int_{0}^{t}  \left(t \wedge \check{\tau}_s - s\right)  \nabla \sigma^{i0}\left(X^{NV}_s\right).~d\vartheta_s \\
& + \frac12  \sum \limits_{k=1}^n \sum \limits_{l=1}^n \int_{0}^{t} \displaystyle \left(t \wedge \check{\tau}_s - s\right) \partial^2_{x_k x_l} \sigma^{i0}\left(X^{NV}_s\right) d\langle \vartheta^k,\vartheta^l \rangle_s\\
&+  \frac12 \sum \limits_{j=1}^d \sum \limits_{k=1}^n \sum \limits_{l=1}^n \int_{0}^{t} \displaystyle \left(t \wedge \check{\tau}_s - s\right) \partial^2_{x_k x_l} \sigma^{i0}\left(X^{NV}_s\right) \sigma^{lj}\left(X^{NV}_s\right) d\langle \vartheta^k,W^j \rangle_s.
\end{split}
\end{equation*}
Now, it is easy to see that \eqref{IPP1_2} is a straightforward consequence of Lemma \ref{Lemme0} and the regularity assumption on the vector fields $\sigma^j$ for $j \in \left\{0,\ldots,d\right\}$.
\end{adem}
We easily obtain the following decomposition of $X^{NV}$
\begin{equation}
\begin{split}
\label{Dec_Xnv}
X^{NV}_t & =\displaystyle \int_{0}^{t}  b\left(X^{NV}_s\right) ds + \sum \limits_{j=1}^d \displaystyle \int_{0}^{t}  \sigma^j\left(X^{NV}_s\right) dW^j_s + \frac{1}{2} \sum \limits_{j=1}^d  \displaystyle \int_{0}^{t} \left(\Delta s - \left(t \wedge \check{\tau}_s - s\right) \right)\left[ \sigma^j, \sigma^0\right ]\left(X^{NV}_s\right) dW^j_s   \\
&+ \frac{1}{4} \displaystyle \int_{0}^{t} \left(\Delta s - \left(t \wedge \check{\tau}_s - s\right) \right)  \left( \partial \sigma^0 \sigma^0 + \sum \limits_{j=1}^d \left( \partial^2\sigma^0 \odot \sigma^j + \partial\sigma^0\partial\sigma^j  \right)\sigma^j  \right)\left(X^{NV}_{\hat{\tau}_s}\right) ds\\
& + \bar{Q}^{2,N}_t + \bar{Q}^{1,N}_t,
\end{split}
\end{equation}
where
\begin{equation}
\begin{split}
\bar{Q}^{2,N}_t &= \frac12 \left(  \displaystyle  \int_{0}^{t} \displaystyle \left(\sigma^0\left(X^{NV}_s\right)  - \sigma^0\left(\bar{X}^0_{s}  \right) \right) ds  - I^0_t \right)  + \frac12 \sum \limits_{j=1}^d \left( \displaystyle \int_{0}^t   \Delta W^j_s \partial \sigma^j \sigma^0\left(X^{NV}_{\hat{\tau}_s}\right)ds -I^j_t \right)  
\end{split}
\end{equation}
is a remainder term with strong order $3/2$ from Proposition \ref{LemmeIPP}. Actually, in \eqref{Dec_Xnv}, the term
\begin{equation}
\label{=0}
\displaystyle \int_{0}^{t} \left(\Delta s - \left(t \wedge \check{\tau}_s - s\right) \right)  \left( \partial \sigma^0 \sigma^0 + \sum \limits_{j=1}^d \left( \partial^2\sigma^0 \odot \sigma^j + \partial\sigma^0\partial\sigma^j  \right)\sigma^j  \right)\left(X^{NV}_{\hat{\tau}_s}\right) ds
\end{equation}
 is null. To lighten up this expression, we denote $F =   \displaystyle \partial \sigma^0 \sigma^0 + \sum \limits_{j=1}^d \left( \partial^2\sigma^0 \odot \sigma^j + \partial\sigma^0\partial\sigma^j  \right)\sigma^j $. Then \eqref{=0} becomes:
\begin{equation}
\begin{split}
\displaystyle \int_{0}^{t} \left(\Delta s - \left(t \wedge \check{\tau}_s - s\right) \right) F\left(X^{NV}_{\hat{\tau}_s}\right) ds&=
\sum \limits_{k=0}^{\lfloor\frac{Nt}{T} \rfloor - 1} F\left(X^{NV}_{t_k}\right) \displaystyle \int_{t_k}^{t_{k+1}} \left(2s - t_k - t_{k+1}\right)  ds\\
& +  F\left(X^{NV}_{\hat{\tau}_t}\right) \displaystyle \int_{\hat{\tau}_t}^{t} \left(2s - \hat{\tau}_t - t\right) ds \\
& = 0.
\end{split}
\end{equation} 	
Therefore \eqref{Dec_Xnv} can be simplified to our final decomposition:
\begin{equation}
\begin{split}
\label{Dec_Xnv_Final}
X^{NV}_t & =\displaystyle \int_{0}^{t}  b\left(X^{NV}_s\right) ds + \sum \limits_{j=1}^d \displaystyle \int_{0}^{t}  \sigma^j\left(X^{NV}_s\right) dW^j_s  + \bar{J}^N_t  + \bar{Q}^{N}_t, 
\end{split}
\end{equation}
where  
\begin{equation}
\bar{J}^N_t =  \frac{1}{2} \sum \limits_{j=1}^d  \displaystyle \int_{0}^{t} \left(\Delta s - \left(t \wedge \check{\tau}_s - s\right) \right)\left[ \sigma^j, \sigma^0\right ]\left(X^{NV}_s\right) dW^j_s, 
\end{equation}
and 
\begin{equation}
\bar{Q}^{N}_t =  \bar{Q}^{2,N}_t +  \bar{Q}^{1,N}_t.
\end{equation}
According to Propositions \ref{Prop} and \ref{LemmeIPP}, $\bar{Q}^{N}$ is a remainder term with strong order $3/2$.
\begin{aprop}
\label{Rest_Q}
Under the assumptions of Theorem \ref{T_EP_CC}: $\forall p \ge 1,\exists C_{5} \in \mathbb{R}_+^*, \forall N \in \mathbb{N}^*,$
\begin{equation}
\label{SCV_cont}
\mathbb{E}\left[ \underset{t\leq T}{\sup} \left\| \bar{Q}^{N}_t \right\|^{2p} \right] \leq C_{5} h^{3p}.
\end{equation}
\end{aprop}
In the following proposition we state the continuous version of Theorem \ref{SC2}.
\begin{aprop}
\label{prop_SCV}
Under the assumptions of Theorem \ref{T_EP_CC}: $\forall p \ge 1,\exists C^{\prime \prime}_{NV} \in \mathbb{R}_+^*, \forall N \in \mathbb{N}^*,$
\begin{equation}
\label{SCV_cont}
\mathbb{E}\left[ \underset{t\leq T}{\sup} \left\|  X_t - X^{NV}_t \right\|^{2p} \right] \leq C^{\prime \prime}_{NV} h^{2p}.
\end{equation}
\end{aprop}
\begin{adem}
Let $p \in [1,+\infty)$, $t \in [0,T]$ and $s\in[0,t]$. Subtracting \eqref{Dec_Xnv_Final} from \eqref{EDS_ITO}, we can evaluate the difference between the exact solution and the scheme: 
\begin{equation}
\begin{split}
X_s - X^{NV}_s & =\displaystyle \int_{0}^{s} \left(b\left(X_u\right) -  b\left(X^{NV}_u\right) \right) du + \sum \limits_{j=1}^d \displaystyle \int_{0}^{s}  \left(\sigma^j\left(X_u\right) - \sigma^j\left(X^{NV}_u\right) \right)dW^j_u  - \bar{J}^N_s  - \bar{Q}^{N}_s. 
\end{split}
\end{equation}
Using a convexity inequality, taking the expectation of the supremum and applying the Burkholder-Davis-Gundy inequality, we get a constant $\alpha_0\in \mathbb{R}_+^*$ independent of $N$, such that 
\begin{equation}
\begin{split}
\mathbb{E}\left[ \underset{s\leq t}{\sup} \left\|  X_s - X^{NV}_s \right\|^{2p} \right] &\leq \alpha_0 \Bigg( \displaystyle \int_{0}^{t} \mathbb{E}\left[ \left\| b\left(X_u\right) -  b\left(X^{NV}_u\right) \right\|^{2p}\right] du + \sum \limits_{j=1}^d \displaystyle \int_{0}^{t} \mathbb{E}\left[ \left\|\sigma^j\left(X_u\right) - \sigma^j\left(X^{NV}_u\right) \right\|^{2p}\right] du\\
&+ \mathbb{E}\left[ \underset{s\leq t}{\sup} \left\| \bar{J}^N_s \right\|^{2p}\right]  + \mathbb{E}\left[ \underset{s\leq t}{\sup} \left\| \bar{Q}^{N}_s \right\|^{2p}\right]\Bigg) . 
\end{split}
\end{equation}
By the Lipschitz assumption, we can rewrite the last inequality as follows 
\begin{equation}
\begin{split}
\label{Gronw}
\mathbb{E}\left[ \underset{s\leq t}{\sup} \left\|  X_s - X^{NV}_s \right\|^{2p} \right] &\leq \alpha_1 \Bigg( \displaystyle \int_{0}^{t} \mathbb{E}\left[ \underset{v\leq u}{\sup} \left\|  X_v - X^{NV}_v \right\|^{2p} \right] du + \mathbb{E}\left[ \underset{s\leq t}{\sup} \left\| \bar{J}^N_s \right\|^{2p}\right]  + \mathbb{E}\left[ \underset{s\leq t}{\sup} \left\| \bar{Q}^{N}_s \right\|^{2p}\right]\Bigg),
\end{split}
\end{equation}
where $\alpha_1\in \mathbb{R}_+^*$ is a constant independent of $N$. On the one hand, applying the Burkholder-Davis-Gundy inequality, we obtain a constant $\beta_0\in \mathbb{R}_+^*$ independent of $N$, such that
\begin{equation}
\mathbb{E}\left[ \underset{s\leq t}{\sup} \left\| \bar{J}^N_s \right\|^{2p}\right]  \leq \beta_0 \sum \limits_{j=1}^d  \displaystyle \int_{0}^{t} \left|\Delta s - \left(t \wedge \check{\tau}_s - s\right) \right|^{2p} \left\| \left[ \sigma^j, \sigma^0\right ]\left(X^{NV}_s\right)\right\|^{2p} du. 
\end{equation}
Since $\left[ \sigma^j, \sigma^0\right ]\in LIP_{loc}^{pgc}\left(\mathbb{R}^n \right)$, \eqref{NV_moment} from Lemma \ref{Lemme0} ensures that
\begin{equation}
\label{GronwJ}
\mathbb{E}\left[ \underset{s\leq t}{\sup} \left\| \bar{J}^N_s \right\|^{2p}\right]  \leq \beta_1 h^{2p}, 
\end{equation}
for some constant $\beta_1\in \mathbb{R}_+^*$ independent of $N$.
On the other hand, from Proposition \ref{Rest_Q}, we have 
\begin{equation}
\label{GronwQ}
\mathbb{E}\left[ \underset{s\leq t}{\sup} \left\| \bar{Q}^{N}_s \right\|^{2p}\right] \leq C_5 h^{3p} 
\end{equation}
Then, combing \eqref{Gronw}, \eqref{GronwJ} and \eqref{GronwQ}, we easily obtain   
\begin{equation}
\begin{split}
\mathbb{E}\left[ \underset{s\leq t}{\sup} \left\|  X_s - X^{NV}_s \right\|^{2p} \right] &\leq \alpha_1 \displaystyle \int_{0}^{t} \mathbb{E}\left[ \underset{v\leq u}{\sup} \left\|  X_v - X^{NV}_v \right\|^{2p} \right] du + \alpha_1\left(\beta_1 + C_5 T^p\right)h^{2p}.
\end{split}
\end{equation}
We conclude thanks to Gr\"{o}nwall's inequality. 
\end{adem}
As a comparison with theorem 4.2 in \cite{NV2}, note that using the adapted interpolation leads to stronger assumptions, which justifies the use of the non-adapted interpolation in \cite{NV2}.

\section{Proof of the stable convergence}
Using \eqref{Dec_Xnv_Final}, the normalized error process $U^N$ can be written as:
\begin{equation}
\begin{split}
\label{Dec_Un}
U^N_t & = N \Bigg(\displaystyle \int_{0}^{t}  \left(b\left(X_s\right) - b\left(X^{NV}_s\right)\right) ds + \sum \limits_{j=1}^d \displaystyle \int_{0}^{t} \left( \sigma^j\left(X_s\right) -  \sigma^j\left(X^{NV}_s\right)\right) dW^j_s\Bigg) + J^N_t  + Q^{N}_t,
\end{split}
\end{equation}
where
\begin{equation}
J^N_t = - N \bar{J}^N_t =  \frac{N}{2} \sum \limits_{j=1}^d  \displaystyle \int_{0}^{t} \left(\Delta s - \left(t \wedge \check{\tau}_s - s\right) \right)\left[ \sigma^0, \sigma^j\right ]\left(X^{NV}_s\right) dW^j_s, 
\end{equation}
and 
\begin{equation}
Q^N =  -N \bar{Q}^N.
\end{equation}
As previously mentioned, by analyzing the stable convergence in law of the source term $J^N$, we prove the stable convergence in law of $U^N$.
In the following, we provide a detailed proof of Theorem \ref{T_EP_CC}. 
\begin{_adem} \textbf{of Theorem \ref{T_EP_CC}:} For the reader's convenience, the proof will go through several steps.
\textbf{Step 1: linearization of $ N  \displaystyle \int_{0}^{t} \displaystyle \left(b\left(X_s\right) - b\left(X^{NV}_s\right)\right) ds  +  N \sum \limits_{j=1}^d \int_{0}^{t} \displaystyle \left( \sigma^j\left(X_s\right) -  \sigma^j\left(X^{NV}_s\right)\right) dW^j_s $. }\\\\
Let $j \in \left\{1,\ldots,d\right\}$ and $i \in \left\{1,\ldots,n\right\}$, by the mean value theorem, we get
\begin{equation*}
\sigma^{ij}\left(X_s\right) -  \sigma^{ij}\left(X^{NV}_s\right) = \nabla \sigma^{ij} \left( \zeta^{ij}_s \right). \left(X_s - X^{NV}_s \right) 
\end{equation*}
where $\zeta^{ij}_s = \alpha^{ij}_s X_s + \left(1 - \alpha^{ij}_s\right) X^{NV}_s$ for some $\alpha^{ij}_s \in [0,1]$.
Using a compact matrix notation, we can write: 
\begin{equation}
 \sigma^{j}\left(X_s\right) - \sigma^{j}\left(X^{NV}_s\right) = \partial \sigma^{j,N}_s \left(X_s - X^{NV}_s \right)
\end{equation} 
where:
\begin{equation}
 \left(\partial \sigma^{j,N}_s\right)_{i,m} = \partial_{x_m} \sigma^{ij}\left(\zeta^{ij}_s\right). 
\end{equation}
In the same way,
\begin{equation*}
b\left(X_s\right) -  b\left(X^{NV}_s\right) = \partial b^N_s \left(X_s - X^{NV}_s \right) 
\end{equation*}
where
\begin{equation}
 \left(\partial b^N_s\right)_{i,m} = \partial_{x_m} b^{i}\left(\zeta^{i0}_s\right) 
\end{equation}
with $\zeta^{i0}_s  = \alpha^{i0}_s X_s + \left(1 - \alpha^{i0}_s\right) X^{NV}_s$ for some $\alpha^{i0}_s \in [0,1]$.
Then, it follows that
\begin{equation*}
\begin{split}
 N  \int_{0}^{t} \displaystyle \left(b\left(X_s\right) - b\left(X^{NV}_s\right)\right) ds  = \int_{0}^{t} \displaystyle \partial b^N_s  U^N_s ds 
\end{split}
\end{equation*}
and
\begin{equation*}
 N \sum \limits_{j=1}^d \int_{0}^{t} \displaystyle \left( \sigma^j\left(X_s\right) -  \sigma^j\left(X^{NV}_s\right)\right) dW^j_s  =  \sum \limits_{j=1}^d \int_{0}^{t} \partial \sigma^{j,N}_s  U^N_s dW^j_s.
\end{equation*}
Then, we can write the error process in three parts as follows:
\begin{equation*}
U^N_t  =  Q^N_t + J^N_t + \left( \int_{0}^{t} \displaystyle \partial b^N_s  U^N_s ds  + \sum \limits_{j=1}^d \int_{0}^{t} \partial \sigma^{j,N}_s  U^N_s dW^j_s\right). 
\end{equation*}

\textbf{Step 2:  stable convergence in law of the source term $J^N$.}\\\\
To study the convergence of the source term $J^N$, we introduce the $d$-dimensional martingale $M^N$ with coordinates
\begin{equation}
M^{j,N}_t = N \int_{0}^{t} \frac12  \left(\Delta s - \left(t \wedge \check{\tau}_s - s\right) \right)  dW^j_s,~~ j \in \left\{1,\ldots,d\right\},
\end{equation}
and $K^N =  \left(K^{j,N} = \left[ \sigma^0, \sigma^j\right ]\left(X^{NV}_s\right)\right)_{1\leq j \leq d}$ with values in $\mathbb{R}^n \otimes \mathbb{R}^d$,
so that
\begin{equation}
J^N_t = \sum \limits_{j=1}^{d} \displaystyle  \int_0^t \displaystyle K_s^{j,N} dM^{j,N}_s.
\end{equation}
\textbf{Step 2.1: stable convergence in law of $M^N$.}\\\\
By virtue of Theorem 2.3 in \cite{NV2}, to study the limit in law of $M^N$, we check $\forall t \in [0,T], \forall j,m,k \in \left\{1,\ldots,d\right\}$
\begin{itemize}
\item the convergence in probability of $\langle M^{j,N},M^{m,N} \rangle_t$, as $N$ goes to infinity,
\item the convergence in probability to $0$ of $\langle M^{j,N},W^k \rangle_t$, as $N$ goes to infinity.  
\end{itemize}
If $j\neq m$ and $j\neq k$  then, obviously,  $\langle M^{j,N},M^{m,N}\rangle_t  =  \langle M^{j,N},W^k \rangle_t = 0.$ Now, if $j=m$, a straightforward calculation gives us:
\begin{equation}
\begin{split}
\langle  M^{j,N},M^{j,N} \rangle_t &=  N^2 \int_{0}^{t} \frac14 \left(\Delta s - \left(t \wedge \check{\tau}_s - s\right) \right)^2 ds \\
&=  \frac14 N^2 \left(\int_{0}^{\hat{\tau}_t } \left(\Delta s -  \left( \check{\tau}_s - s\right) \right)^2 ds + \int_{\hat{\tau}_t }^{t}  \left( s - \hat{\tau}_t - \left(t - s\right) \right)^2 ds\right)\\
&=  \frac1{12} N^2 \left(\lfloor\frac{Nt}{T} \rfloor  \frac{T^3}{N^3}+  \left(t-\hat{\tau}_t \right)^3\right) \underset{N \to +\infty}{\longrightarrow} \frac1{12} t T^2.
\end{split}
\end{equation}
If $j= k$ 
\begin{equation}
\begin{split}
\langle M^{j,N},W^k \rangle_t &= N \int_{0}^{t} \frac12 \left(\Delta s - \left(t \wedge \check{\tau}_s - s\right) \right) ds = 0.
\end{split}
\end{equation}
Applying  Theorem 2.3 in \cite{NV2} we conclude that $\displaystyle \frac{2\sqrt{3}}{T}M^N$ stably converges in law to a standard $d$-dimensional Brownian motion $\tilde{B}$, independent of $W$.\\\\
\textbf{Step 2.2: convergence in probability of $K^N $.}\\\\
Now, it remains to prove the convergence in probability of $K^N$. From Proposition \ref{prop_SCV}, together with the continuity assumption on $\left[\sigma^0,\sigma^j\right]$, $j \left\{1,\ldots,d\right\}$, we get the following convergence in probability
\begin{equation}
\underset{t\leq T}{\sup} \left\| K_t^{j,N} - \left[\sigma^0,\sigma^j\right]\left(X_{t}\right) \right\| \overset{\mathbb{P}}{\underset{N \to +\infty}{\longrightarrow}} 0.
\end{equation}
\textbf{Step 2.3: conclusion of the step 2.}\\\\
According to Proposition 2.2 in \cite{NV2}, we have the following convergence:
\begin{equation}
\left(K^N, \frac{2\sqrt{3}}{T}M^N\right) \overset{stably}{\underset{N \to +\infty}{\Longrightarrow}} \left(\left(\left[\sigma^0,\sigma^j\right]\left(X\right)\right)_{j \in \left\{1,\ldots,d\right\}},\tilde{B}\right).
\end{equation}
The convergence of $\left<M^N\right>_T$ ensures its tightness. Then Proposition 2.4 in \cite{NV2} leads us to:
\begin{equation}
\label{Step3}
\left(K^N, \frac{2\sqrt{3}}{T}M^N, J^N\right) \overset{stably}{\underset{N \to +\infty}{\Longrightarrow}} \left(\left(\left[\sigma^0,\sigma^j\right]\left(X\right)\right)_{j \in \left\{1,\ldots,d\right\}},\tilde{B},\left(\frac{T}{2\sqrt{3}} \sum \limits_{j=1}^d  \int_0^t \displaystyle \left[\sigma^0,\sigma^j\right]\left(X_s\right) d\tilde{B}^{j}_s \right)_{t\in [0,T]} \right).
\end{equation}
\textbf{Step 3: convergence of $Q^N$.}\\\\
We easily get the following convergence in $L^2$ from Proposition \ref{Rest_Q}:
\begin{equation}
\label{Step4}
\underset{t\leq T}{\sup} \left\| Q_t^{N} \right\| \overset{L^2}{\underset{N \to +\infty}{\longrightarrow}} 0.
\end{equation}
\\\\

\textbf{Step 4: stable convergence in law of $U^N$.}\\\\
We recall that $U^{N}_t = Q^N_t + J^N_t + \left(\displaystyle \int_0^t  \partial b_s^NU^N_s ds + \sum \limits_{j=1}^d \displaystyle \int_0^t \displaystyle \partial \sigma^{j,N}_s U^N_s dW_s^j \right)$.
Thanks to \eqref{Step3} and \eqref{Step4}, we conclude using Theorem 2.5 in \cite{NV2} since the continuity of $\partial b$ and $\partial \sigma^j, j \in  \left\{1,\ldots,d\right\},$ together with  Proposition \ref{prop_SCV}, ensure that 
\begin{equation}
\underset{t\leq T}{\sup} \left\| \partial b^N_t - \partial b\left(X_t\right) \right\| \overset{\mathbb{P}}{\underset{N \to +\infty}{\longrightarrow}} 0,
\end{equation}
and $\forall j \in \left\{1,\ldots,d\right\}$,
\begin{equation}
\underset{s\leq T}{\sup} \left\| \partial \sigma^{j,N}_t - \partial \sigma^{j}\left(X_t\right) \right\| \overset{\mathbb{P}}{\underset{N \to +\infty}{\longrightarrow}} 0.
\end{equation}
\end{_adem}
\section{Appendix}
This section is devoted to the proof of Propositions \ref{Prop_App_1}, \ref{Prop_App_2} and \ref{Prop_Q}. Before proving Proposition \ref{Prop_App_1}, in the following lemma, we give intermediate estimations of the process $Y^{d+1}$.
\begin{alem}
\label{Lemme2_A}
Assume that $\forall j \in \left\{0,\ldots,d\right\}, \sigma^j \in \mathcal{C}^{1}\left(\mathbb{R}^n,\mathbb{R}^n\right)$  with bounded first order derivatives. 
 Then, $\forall p \ge 1, \exists C_6 \in \mathbb{R}_+^*, \forall N \in \mathbb{N}^*, \forall t \in [0,T], \forall s \in [\hat{\tau}_{t},t],$
\begin{equation}
\begin{split}
\label{Acc_Yd+1_bis}
\mathbb{E}\left[ \left\|  Y_{t,s}^{d+1}  - Id_n \right\|^{2p}\right] & \leq   C_6 \left(\Delta t\right)^{2p},
\end{split}
\end{equation}
\begin{equation}
\begin{split}
\label{Acc_Yd+1}
\mathbb{E}\left[ \left\|  Y_{t,t}^{d+1} \left(Y^{d+1}_{t,s}\right)^{-1} - Id_n \right\|^{2p}\right] & \leq   C_6 \left(\Delta t\right)^{2p}.
\end{split}
\end{equation}
\end{alem}
\begin{adem}
Let $t \in [0,T]$ and $s \in [\hat{\tau}_{t},t]$, we recall that $Y^{d+1}_{t,s}$ is the solution, at time $s$, of the following ODE:
\begin{equation}
\left\{
    \begin{array}{ll}
\frac{d}{du} \Xi_{t,\hat{\tau}_{t},u} =  \frac12 \partial \sigma^0 \left( \bar{X}^{d+1}_{t,u} \right)   \Xi_{t,\hat{\tau}_{t},u}, u \in [\hat{\tau}_{t},t]\\
 \Xi_{t,\hat{\tau}_{t},\hat{\tau}_{t}} = Id_n,
\label{EDOYd+1}
\end{array}
\right.
\end{equation}
and $Y_{t,t}^{d+1} \left(Y^{d+1}_{t,s}\right)^{-1}$ is the solution, at time $t$, of the same ODE, but with the initial condition:
\begin{equation*}
\left\{
    \begin{array}{ll}
\frac{d}{du} \Xi_{t,s,u} =  \frac12 \partial \sigma^0 \left( \bar{X}^{d+1}_{t,u} \right)   \Xi_{t,\hat{\tau}_{t},u}, u \in [s,t]\\
 \Xi_{t,s,s} = Id_n.
\end{array}
\right.
\end{equation*}
Therefore, we deduce \eqref{Acc_Yd+1} similarly to \eqref{Acc_Yd+1_bis}.
Since $\partial \sigma^0$ is bounded, there exists a constant $c \in \mathbb{R}_+^*$, independent of $N$, $t$ and $s$, such that
\begin{equation}
\begin{split}
\label{Estim_EDO_Yd+1}
\left\| Y^{d+1}_{t,s} \right\|^{2p}   \leq \exp\left(c \left(s -\hat{\tau}_{t}\right)\right) \leq \exp\left(c T\right),
\end{split} 
\end{equation}
Writing $Y^{d+1}_{t,s}$ in integral form, we have
\begin{equation}
\begin{split}
\label{EDO_1}
Y^{d+1}_{t,s} &=    Id_n + \frac12 \int_{\hat{\tau}_{t}}^{s} \displaystyle \partial \sigma^0 \left( \bar{X}^{d+1}_{t,u} \right)  Y^{d+1}_{t,u}du.
\end{split} 
\end{equation} 
Hence, using a convexity inequality,
\begin{equation}
\begin{split}
\mathbb{E}\left[ \left\|  Y_{t,s}^{d+1}  - Id_n \right\|^{2p}\right] & \leq  \frac1{2^{2p}} \exp\left(c T\right) \left(\Delta t\right)^{2p-1} \int_{\hat{\tau}_{t}}^{s} \displaystyle \left\| \partial \sigma^0 \left( \bar{X}^{d+1}_{t,u} \right) \right\|^{2p} du.
\end{split}
\end{equation}
Since $\partial \sigma^0$ is bounded we easily get \eqref{Acc_Yd+1_bis}. 
\end{adem}
\begin{_adem} \textbf{of Proposition \ref{Prop_App_1}:}\\\\
\textbf{Proof of \eqref{approx_1}:}\\\\ 
To prove \eqref{approx_1}, we only need to assume that
\begin{itemize}
\item $\forall j \in \left\{0,\ldots,d\right\}, \sigma^j \in \mathcal{C}^{1}\left(\mathbb{R}^n,\mathbb{R}^n\right)$  with bounded first order derivatives,
\item $\sum \limits_{j=1}^d \partial \sigma^j \sigma^j$ is a Lipschitz continuous function,
\item $\partial \sigma^0 \in  LIP_{loc}^{pgc}\left(\mathbb{R}^n \otimes \mathbb{R}^n\right)$.
\end{itemize}
Let $t \in [0,T]$, $s \in [\hat{\tau}_{t},t]$, $p \ge 1$, and
\begin{equation}
\label{Theta_d+1}
\Theta^{d+1}_{t,s} = Y^{d+1}_{t,s} -  Id_n - \frac{1}{2} \Delta s  \partial \sigma^0\left(X^{NV}_t\right).
\end{equation} 
Adding and subtracting $\displaystyle  \frac12 \int_{\hat{\tau}_{t}}^{s} \displaystyle \partial \sigma^0 \left( X^{NV}_t \right)  Y^{d+1}_{t,u}du$, we get
\begin{equation}
\begin{split}
\Theta^{d+1}_{t,s} &= \frac12 \int_{\hat{\tau}_{s}}^{s} \displaystyle \left( \partial \sigma^0 \left( \bar{X}^{d+1}_{t,u} \right) - \partial \sigma^0\left(X^{NV}_t\right) \right) Y^{d+1}_{t,u}du + \frac12 \int_{\hat{\tau}_{s}}^{s} \displaystyle  \partial \sigma^0\left(X^{NV}_t\right)  \left(Y^{d+1}_{t,u} - Id_n \right)du. 
\end{split}
\end{equation} 
Combining a convexity inequality and the Cauchy-Schwarz inequality, we obtain
\begin{equation}
\begin{split}
\mathbb{E}\left[ \left\|  \Theta^{d+1}_{t,s}\right\|^{2p}\right]  &\leq   \frac{2^{2p-1}}{2^{2p}}  \left(\Delta t\right)^{2p-1}  \Bigg( \int_{\hat{\tau}_{s}}^{s} \displaystyle \left(\mathbb{E}\left[ \left\| \partial \sigma^0 \left( \bar{X}^{d+1}_{t,u} \right) - \partial \sigma^0\left(X^{NV}_t\right) \right\|^{4p} \right] \mathbb{E}\left[ \left\| Y^{d+1}_{t,u} \right\|^{4p} \right]  \right)^{\frac12} du \\
& + \int_{\hat{\tau}_{s}}^{s} \displaystyle \left(\mathbb{E}\left[ \left\| \partial \sigma^0\left(X^{NV}_t\right) \right\|^{4p} \right] \mathbb{E}\left[ \left\| Y^{d+1}_{t,u} - Id_n \right\|^{4p}\right]  \right)^{\frac12} du \Bigg). 
\end{split}
\end{equation} 
Since $\partial \sigma^0$ is bounded and locally Lipschitz with polynomially growing Lipschitz constant, there exists a constant $c \in \mathbb{R}_+^*$ and $q \in \mathbb{N}$, independent of $N$, $t$ and $s$, such that 
\begin{equation}
\begin{split}
\mathbb{E}\left[ \left\|  \Theta^{d+1}_{t,s}\right\|^{2p}\right]  &\leq   c  \left(\Delta t\right)^{2p-1} \Bigg(  \int_{\hat{\tau}_{s}}^{s} \displaystyle \left(\mathbb{E}\left[ \left( 1 + \left\| \bar{X}^{d+1}_{t,u} \right\|^{q} \vee \left\|  X^{NV}_t \right\|^{q}\right)  \left\| \bar{X}^{d+1}_{t,u} - X^{NV}_t \right\|^{4p} \right] \mathbb{E}\left[ \left\| Y^{d+1}_{t,u} \right\|^{4p} \right]  \right)^{\frac12} du\\
&  + \int_{\hat{\tau}_{s}}^{s} \displaystyle \left(\mathbb{E}\left[ \left\| Y^{d+1}_{t,u} - Id_n \right\|^{4p}\right]  \right)^{\frac12} du \Bigg). 
\end{split}
\end{equation}
Applying once again the Cauchy-Schwarz inequality, and using \eqref{NV_moment} and \eqref{Acc_d+1} from Lemma \ref{Lemme0}, and \eqref{Estim_EDO_Yd+1} from the above proof of Lemma \ref{Lemme2_A}, we get a constant $\alpha_1 \in \mathbb{R}_+^*$, independent of $N$, $t$ and $s$, such that
\begin{equation}
\begin{split}
\mathbb{E}\left[ \left\|  \Theta^{d+1}_{t,s}\right\|^{2p}\right]  &\leq   \alpha_1   \left(\Delta t\right)^{2p-1}  \int_{\hat{\tau}_{s}}^{s} \displaystyle \left(  \left(\mathbb{E}\left[   \left\| \bar{X}^{d+1}_{t,u} - X^{NV}_t \right\|^{8p} \right]  \right)^{\frac14} + \left(\mathbb{E}\left[ \left\| Y^{d+1}_{t,u} - Id_n \right\|^{8p}\right]  \right)^{\frac14}\right) du . 
\end{split}
\end{equation}
We conclude using \eqref{Acc_d+1} from Lemma \ref{Lemme0}, and \eqref{Acc_Yd+1_bis}  from Lemma \ref{Lemme2_A}.\\

\textbf{Proof of \eqref{approx_4}:}\\\\ 
To prove \eqref{approx_4}, we only need to assume that
\begin{itemize}
\item $\forall j \in \left\{1,\ldots,d\right\}, \sigma^j \in \mathcal{C}^{1}\left(\mathbb{R}^n,\mathbb{R}^n\right)$  with bounded first order derivatives,
\item $\sum \limits_{j=1}^d \partial \sigma^j \sigma^j$ is a Lipschitz continuous function,
\item $\sigma^0 \in \mathcal{C}^{2}\left(\mathbb{R}^n,\mathbb{R}^n\right)$  with bounded first order derivatives and $\partial^2 \sigma^0 \in  LIP_{loc}^{pgc}\left(\mathbb{R}^n \otimes \mathbb{R}^n\right)$.
\end{itemize}

Let $t \in [0,T]$, $j \in \left\{1,\ldots,d\right\}$, $p \ge 1$, and
\begin{equation}
 \theta^{d+1,j}_t = Z^{d+1,j}_{t} - \frac{1}{2} \Delta t \left( \partial^2\sigma^0\odot \sigma^j\right) \sigma^j\left(X^{NV}_{\hat{\tau}_{t}}\right) ,
\end{equation} 
and we recall that 
\begin{equation*}
Z^{d+1,j}_t = \frac{1}{2}\int_{\hat{\tau}_{t}}^{t} \displaystyle  Y_{t,t}^{d+1} \left(Y^{d+1}_{t,s}\right)^{-1} \displaystyle \left( \partial^2 \sigma^0\left(\bar{X}^{d+1}_{t,s}\right) \odot \left(Y^{d+1}_{t,s} \sigma^j\left(\bar{X}_{t,t}\right)\right)\right) Y^{d+1}_{t,s} \sigma^j\left(\bar{X}_{t,t}\right)\ ds.
\end{equation*}
Therefore,
\begin{equation}
\begin{split}
\theta^{d+1,j}_t & = \frac{1}{2}\int_{\hat{\tau}_{t}}^{t} \displaystyle \left(  Y_{t,t}^{d+1} \left(Y^{d+1}_{t,s}\right)^{-1}\displaystyle \left( \partial^2 \sigma^0\left(\bar{X}_{t,s}\right) \odot Y^{d+1}_{t,s} \sigma^j\left(\bar{X}_{t,t}\right)\right) Y^{d+1}_{t,s} \sigma^j\left(\bar{X}_{t,t}\right) -  \left( \partial^2\sigma^0\odot \sigma^j\right) \sigma^j\left(X^{NV}_{\hat{\tau}_{t}}\right) \right) ds.   
 \end{split}
\end{equation} 
Adding and subtracting some appropriate terms, we obtain:
\begin{equation}
\begin{split}
\theta^{d+1,j}_t &= \frac{1}{2}\int_{\hat{\tau}_{t}}^{t} \displaystyle   Y^{d+1}_{t,t} \left( Y_{t,s}^{d+1}\right)^{-1}\displaystyle \left( \partial^2 \sigma^0\left(\bar{X}_{t,s}\right) \odot Y^{d+1}_{t,s} \sigma^j\left(\bar{X}_{t,t}\right)\right) \left(Y^{d+1}_{t,s} - Id_n\right) \sigma^j\left(\bar{X}_{t,t}\right) ds  \\
& + \displaystyle \frac{1}{2}\displaystyle  \int_{\hat{\tau}_{t}}^{t} \displaystyle \left(   Y_{t,t}^{d+1}\left(Y^{d+1}_{t,s}\right)^{-1} - Id_n\right) \displaystyle \left( \partial^2 \sigma^0\left(\bar{X}_{t,s}\right) \odot Y^{d+1}_{t,s} \sigma^j\left(\bar{X}_{t,t}\right)\right)  \sigma^j\left(\bar{X}_{t,s}\right) ds\\
&  + \frac{1}{2}\int_{\hat{\tau}_{t}}^{t} \displaystyle   \displaystyle \left(\partial^2 \sigma^0\left(\bar{X}_{t,s}\right) \odot \left(Y^{d+1}_{t,s} - Id_n\right) \sigma^j\left(\bar{X}_{t,t}\right)\right)\sigma^j\left(\bar{X}_{t,t}\right)  ds \\
& + \frac{1}{2}\int_{\hat{\tau}_{t}}^{t} \displaystyle   \displaystyle \left( \left(\partial^2\sigma^0\left(\bar{X}_{t,s}\right) - \partial^2\sigma^0\left(X^{NV}_t\right) \right) \odot \sigma^j\left(\bar{X}_{t,t}\right)\right) \sigma^j\left(\bar{X}_{t,t}\right) ds \\
& + \frac{1}{2}\int_{\hat{\tau}_{t}}^{t} \displaystyle   \displaystyle \left(  \partial^2\sigma^0\left(X^{NV}_t\right)  \odot \left(\sigma^j\left(\bar{X}_{t,t}\right) -\sigma^j\left(X^{NV}_{t}\right) \right)\right) \sigma^j\left(\bar{X}_{t,t}\right) ds \\
& + \frac{1}{2}\int_{\hat{\tau}_{t}}^{t} \displaystyle   \displaystyle \left(  \partial^2\sigma^0\left(X^{NV}_t\right)  \odot \sigma^j\left(X^{NV}_t\right)\right) \left(\sigma^j\left(\bar{X}_{t,t}\right) -\sigma^j\left(X^{NV}_{t}\right) \right) ds \\
& + \frac{1}{2}\int_{\hat{\tau}_{t}}^{t} \displaystyle   \displaystyle \left( \partial^2\sigma^0\odot \sigma^j\right) \sigma^j\left(X^{NV}_t\right) - \left( \partial^2\sigma^0\odot \sigma^j\right) \sigma^j\left(X^{NV}_{\hat{\tau}_{t}}\right) ds.
 \end{split}
\end{equation} 
Note that since $\partial^2 \sigma^0 \in LIP_{loc}^{pgc}\left(\mathbb{R}^n \otimes \mathbb{R}^n \otimes \mathbb{R}^n\right)$ and $\sigma^j$ is Lipschitz continuous, then $\left( \partial^2\sigma^0\odot \sigma^j\right) \sigma^j \in LIP_{loc}^{pgc}\left(\mathbb{R}^n\right)$. 
We easily get the desired result by combing a convexity inequality, the Cauchy-Schwarz inequality, \eqref{Acc_Yd+1_bis}, \eqref{Acc_Yd+1} and \eqref{Estim_EDO_Yd+1} from Lemma \ref{Lemme2_A}, \eqref{NV_moment} \eqref{Acc_NV} and \eqref{Acc_} from Lemma \ref{Lemme0}

\textbf{Proof of \eqref{approx_2}:}\\\\ 
To prove \eqref{approx_2}, we assume 
\begin{itemize}
\item $\forall j \in \left\{0,\ldots,d\right\}, \sigma^j \in \mathcal{C}^{1}\left(\mathbb{R}^n,\mathbb{R}^n\right)$ with bounded first order derivatives,
\item $\sum \limits_{j=1}^d \partial \sigma^j \sigma^j$ is a Lipschitz continuous function,
\item $\forall j \in \left\{1,\ldots,d\right\}, \partial \sigma^j \in  LIP_{loc}^{pgc}\left(\mathbb{R}^n \otimes \mathbb{R}^n\right)$,
\end{itemize}
and that the commutativity condition \eqref{Comm_condition} holds.
Let $t \in [0,T]$, $j \in \left\{1,\ldots,d\right\}$, $i \in \left\{1,\ldots,n\right\}$, $p \ge 1$, and
\begin{equation}
 \theta^j_t = \sigma^j\left(\bar{X}_{t,t}\right)  - \sigma^j\left(X^{NV}_{t}\right) + \frac{1}{2} \Delta t \partial \sigma^j\sigma^0 \left(X^{NV}_{t} \right).
\end{equation}
The $i$-th component of $\theta^{j}$ is given by:
\begin{equation*}
 \theta^{ij}_t  = \sigma^{ij}\left(\bar{X}_{t,t}\right) - \sigma^{ij}\left(X^{NV}_{t}\right) + \frac{1}{2} \Delta t \nabla \sigma^{ij}\left(X^{NV}_{t} \right).~\sigma^{0} \left(X^{NV}_{t} \right). 
\end{equation*}
By the mean value theorem
\begin{equation}
 \sigma^{ij}\left(\bar{X}_{t,t}\right)  =  \sigma^{ij}\left(X^{NV}_{t}\right) + \nabla \sigma^{ij}\left(X^{NV}_{t} \right). \left(\bar{X}_{t,t} - X^{NV}_{t}\right) + \left( \nabla \sigma^{ij}\left(\xi^{ij}_{t} \right) - \nabla \sigma^{ij}\left(X^{NV}_{t} \right) \right). \left(\bar{X}_{t,t} - X^{NV}_{t}\right) 
\end{equation}
where $\xi^{ij}_{t} = X^{NV}_t + \alpha^{ij}_t \left(\bar{X}_{t,t} - X^{NV}_t\right)$ for some $\alpha^{ij}_t \in [0,1]$.
Then, it follows that
\begin{equation*}
 \theta^{ij}_t  = \nabla \sigma^{ij}\left(X^{NV}_{t} \right). \left(\bar{X}_{t,t} - X^{NV}_{t}\right) + \left( \nabla \sigma^{ij}\left(\xi^{ij}_{t} \right) - \nabla \sigma^{ij}\left(X^{NV}_{t} \right) \right). \left(\bar{X}_{t,t} - X^{NV}_{t}\right)   + \frac{1}{2} \Delta t \nabla \sigma^{ij}\left(X^{NV}_{t} \right).~\sigma^{0} \left(X^{NV}_{t} \right). 
\end{equation*}
Moreover, since
\begin{equation*}
X^{NV}_t = \exp\left(\frac{1}{2}\Delta t \sigma^0 \right) \bar{X}_{t,t},
\end{equation*}
we have that
\begin{equation*}
\bar{X}_{t,t} - X^{NV}_{t} =   - \frac{1}{2} \int_{\hat{\tau}_{t}}^{t} \displaystyle \sigma^0\left(\bar{X}^{d+1}_{t,s}\right) ds.
\end{equation*}
Therefore, we obtain
\begin{equation}
\begin{split}
 \theta^{ij}_t  &=  \frac{1}{2} \int_{\hat{\tau}_{t}}^{t} \nabla \sigma^{ij}\left(X^{NV}_{t} \right).\left(\sigma^{0} \left(X^{NV}_{t} \right) - \sigma^0\left(\bar{X}^{d+1}_{t,s}\right) \right) ds +\left( \nabla \sigma^{ij}\left(\xi^{ij}_{t} \right) - \nabla \sigma^{ij}\left(X^{NV}_{t} \right) \right). \left(\bar{X}_{t,t} - X^{NV}_{t}\right)   . 
\end{split}
\end{equation}
Using a convexity inequality
\begin{equation}
\begin{split}
\mathbb{E}\left[ \left| \theta^{ij}_t \right|^{2p}\right] & \leq   2^{2p-1}  \Bigg( \frac1{2^{2p}}\left(\Delta t\right)^{2p-1}  \int_{\hat{\tau}_{t}}^{t} \mathbb{E}\left[ \left| \nabla \sigma^{ij}\left(X^{NV}_{t} \right).\left(\sigma^{0} \left(X^{NV}_{t} \right) - \sigma^0\left(\bar{X}^{d+1}_{t,s}\right) \right) \right|^{2p}\right] ds\\
& +\displaystyle \mathbb{E}\left[ \left| \left( \nabla \sigma^{ij}\left(\xi^{ij}_{t} \right) - \nabla \sigma^{ij}\left(X^{NV}_{t} \right) \right). \left(\bar{X}_{t,t} - X^{NV}_{t}\right) \right|^{2p}\right]  \Bigg). 
\end{split}
\end{equation}
 Applying the Cauchy-Schwarz inequality, we get
 \begin{equation}
\begin{split}
\mathbb{E}\left[ \left| \theta^{ij}_t \right|^{2p}\right] & \leq 2^{2p-1}  \Bigg( \frac1{2^{2p}}\left(\Delta t\right)^{2p-1}  \int_{\hat{\tau}_{t}}^{t}\left( \mathbb{E}\left[ \left\| \nabla \sigma^{ij}\left(X^{NV}_{t} \right)\right\|^{4p} \right] \mathbb{E}\left[ \left\| \sigma^{0} \left(X^{NV}_{t} \right) - \sigma^0\left(\bar{X}^{d+1}_{t,s}\right)  \right\|^{4p}\right]\right)^{\frac12} ds\\
& +\displaystyle \left( \mathbb{E}\left[ \left\| \nabla \sigma^{ij}\left(\xi^{ij}_{t} \right) - \nabla \sigma^{ij}\left(X^{NV}_{t} \right)\right\|^{4p} \right] \mathbb{E}\left[ \left\| \bar{X}_{t,t} - X^{NV}_{t} \right\|^{4p}\right]\right)^{\frac12}  \Bigg). 
\end{split}
\end{equation}
Since $\sigma^0$ is Lipschitz continuous and $\partial \sigma^j$ is locally Lipschitz with polynomially growing Lipschitz constant, using \eqref{NV_moment} from Lemma \ref{Lemme0} , we easily get a constant $\alpha_2 \in \mathbb{R}_+^*$ independent of $N$ and $t$ such that
\begin{equation}
\begin{split}
\mathbb{E}\left[ \left| \theta^{ij}_t \right|^{2p}\right] & \leq   \alpha_2 \Bigg( \left(\Delta t\right)^{2p-1} \int_{\hat{\tau}_{t}}^{t} \displaystyle \left( \mathbb{E}\left[\left\| X^{NV}_{t}  - \bar{X}^{d+1}_{t,s} \right\|^{4p} \right]\right)^{\frac12} ds   +  \mathbb{E}\left[ \left\|  \bar{X}_{t,t}   - X^{NV}_{t}   \right\|^{4p}\right] \Bigg). 
\end{split}
\end{equation}
Applying \eqref{Acc_d+1} and \eqref{Acc_t} from Lemma \ref{Lemme0} , we obtain
\begin{equation}
\begin{split}
\label{sim}
\mathbb{E}\left[ \left| \theta^{ij}_t \right|^{2p}\right] & \leq   \alpha_2 \left( \sqrt{C_0} + C_0  \right) \left(\Delta t\right)^{4p}.  
\end{split}
\end{equation}

\textbf{Proof of \eqref{approx_3}:}\\\\ 
The proof of \eqref{approx_3} is very similar to \eqref{approx_2}. However, we need to assume further that $\forall j \in \left\{1,\ldots,d\right\}, \sigma^j  \in \mathcal{C}^{2}\left(\mathbb{R}^n,\mathbb{R}^n\right)$ and $\partial^2 \sigma^j \in  LIP_{loc}^{pgc}\left(\mathbb{R}^n \otimes \mathbb{R}^n \otimes \mathbb{R}^n\right)$.
Let $t \in [0,T]$, $j \in \left\{1,\ldots,d\right\}$, $p \ge 1$, and
\begin{equation}
\bar{\theta}^j_t = \partial \sigma^j \sigma^j\left(\bar{X}_{t,t}\right)  - \partial \sigma^j \sigma^j\left(X^{NV}_{t}\right) + \frac{1}{2} \Delta t F_{j} \left(X^{NV}_{\hat{\tau}_t} \right).
\end{equation}
where $ F_{j} = \left(\partial^2 \sigma^j \odot \sigma^j + \left(\partial \sigma^j\right)^2\right)\sigma^0$.
Since $\partial^2 \sigma^j \in LIP_{loc}^{pgc}\left(\mathbb{R}^n \otimes \mathbb{R}^n \otimes \mathbb{R}^n\right)$, $ \partial \sigma^j \in LIP_{loc}^{pgc}\left(\mathbb{R}^n \otimes \mathbb{R}^n\right)$ and $\sigma^0$ is  Lipschitz continuous, then $F_{j} \in LIP_{loc}^{pgc}\left(\mathbb{R}^n\right)$. Hence, similarly to \eqref{sim}, there exists $\bar{\alpha}_2 \in \mathbb{R}_+^*$, independent of $N$ and $t$, such that
\begin{equation}
\begin{split}
\label{idem}
\mathbb{E}\left[ \left\| \partial \sigma^j \sigma^j\left(\bar{X}_{t,t}\right)  - \partial \sigma^j \sigma^j\left(X^{NV}_{t}\right) + \frac{1}{2} \Delta t F_{j} \left(X^{NV}_{t} \right)  \right\|^{2p}\right] & \leq   \bar{\alpha}_2 \left(\Delta t\right)^{4p}.
\end{split}
\end{equation}
Then, using a convexity inequality
\begin{equation}
\begin{split}
\mathbb{E}\left[ \left\| \bar{\theta}^{j}_t \right\|^{2p}\right] &\leq 2^{2p-1} \Bigg( \mathbb{E}\left[ \left\| \partial \sigma^j \sigma^j\left(\bar{X}_{t,t}\right)  - \partial \sigma^j \sigma^j\left(X^{NV}_{t}\right) + \frac{1}{2} \Delta t F_{j} \left(X^{NV}_{t} \right) \right\|^{2p}\right] \\
& + \mathbb{E}\left[ \left\| \frac{1}{2} \Delta t \left( F_{j} \left(X^{NV}_{\hat{\tau}_t} \right)  - F_{j} \left(X^{NV}_{t} \right) \right) \right\|^{2p}\right] \Bigg)\\ 
&\leq 2^{2p-1} \Bigg( \bar{\alpha}_2 \left(\Delta t\right)^{4p} +   \frac{1}{2^{2p}} \left(\Delta t\right)^{2p} \mathbb{E}\left[ \left\| F_{j} \left(X^{NV}_{\hat{\tau}_t} \right)  - F_{j} \left(X^{NV}_{t} \right) \right\|^{2p}\right] \Bigg)
\end{split}
\end{equation}
Since $F_j$ is locally Lipschitz with polynomially growing Lipschitz constant, we conclude using \eqref{Acc_NV} from Lemma \ref{Lemme0}.
\end{_adem}
Before proving Proposition \ref{Prop_App_2}, we introduce some intermediate processes. We define for $j \in \left\{1,\ldots,d\right\}$, $t \in [0,T]$ and $s \in [\hat{\tau}_t,t]$
\begin{equation}
\bar{X}^j_{t,s} =  h_j\left(\frac{\Delta t}{2},\Delta W^1_t,\ldots,\Delta W^{j-1}_t,\Delta W^j_s;X^{NV}_{\hat{\tau}_{t}}\right).
\end{equation}
The next lemma, which is similar to Lemma \ref{Lemme0}, compares the Ninomiya-Victoir scheme to the intermediate process $\bar{X}^j,  j \in \left\{1,\ldots,d\right\}$.
We omit its proof.
\begin{alem}
\label{Lemme0_A}
Assume that 
\begin{itemize}
\item $\sigma^0$ is a Lipschitz continuous function,
\item $\forall j \in \left\{1,\ldots,d\right\}, \sigma^j \in \mathcal{C}^{1}\left(\mathbb{R}^n,\mathbb{R}^n\right)$  with bounded first order derivatives,
\item $\forall j \in \left\{1,\ldots,d\right\}, \partial \sigma^j \sigma^j$ is a Lipschitz continuous function.
\end{itemize}
Then, $\forall p \ge 1, \exists C_7 \in \mathbb{R}_+^*, \forall N \in \mathbb{N}^*, \forall t \in [0,T], \forall s, s_1, s_2 \in [\hat{\tau}_t,t], \forall j,m \in \left\{1,\ldots,d\right\},$
\begin{equation}
\begin{split}
\label{Acc_j_tau}
\mathbb{E}\left[ \left\|  X^{NV}_{\hat{\tau}_t} - \bar{X}^j_{t,\hat{\tau}_t} \right\|^{2p}\right] & \leq   C_7 \left(\Delta t\right)^{p},
\end{split}
\end{equation}
\begin{equation}
\begin{split}
\label{Acc_j}
\mathbb{E}\left[ \left\| X^{NV}_{t}  - \bar{X}^j_{t,s}\right\|^{2p}\right] & \leq   C_7 \left(\Delta t\right)^{p}.
\end{split}
\end{equation}
\begin{equation}
\label{rem2_1}
\mathbb{E}\left[ \left\| \bar{X}^j_{t,s_1} - \bar{X}^{m}_{t,s_2} \right\|^{2p}\right]  \leq   C_7 \left(\Delta t\right)^{p},
\end{equation}
\end{alem}
In order to derive the estimation \eqref{approx_5} from Proposition \ref{Prop_App_2}, we also need the following lemma, which gives several approximations of the processes $Y^j$, for $j \in \left\{1,\ldots,d\right\}$. 
\begin{alem}
\label{Lemme4_A}
Assume that 
\begin{itemize}
\item $\sigma^0$ is a Lipschitz continuous function,
\item $\forall j \in \left\{1,\ldots,d\right\}, \sigma^j \in \mathcal{C}^{2}\left(\mathbb{R}^n,\mathbb{R}^n\right)$  with bounded first order derivatives, $\partial \sigma^j \in LIP_{loc}^{pgc}\left(\mathbb{R}^n \otimes \mathbb{R}^n\right)$ and $\partial^2 \sigma^j  \in LIP_{loc}^{pgc}\left(\mathbb{R}^n \otimes \mathbb{R}^n \otimes \mathbb{R}^n\right)$,
\item $\forall j \in \left\{1,\ldots,d\right\}, \partial \sigma^j \sigma^j$ is a Lipschitz continuous function.
\end{itemize}
 Then, $\forall p \ge 1, \exists C_8 \in \mathbb{R}_+^*, \forall t \in [0,T], \forall s \in [\hat{\tau}_{t},t], \forall N \in \mathbb{N}^*$
\begin{equation}
\begin{split}
\label{Acc_Yj}
\mathbb{E}\left[ \left\|  Y_{t,s}^{j}  - Id_n \right\|^{2p}\right] & \leq   C_8 \left(\Delta t\right)^{p},
\end{split}
\end{equation}
\begin{equation}
\begin{split}
\label{approx_Yj}
\mathbb{E}\left[ \left\| Y^{j}_{t,s}  - Id_n -  \Delta W^j_s  \partial \sigma^j\left(\bar{X}^{j}_{t,\hat{\tau}_{t}}\right) \right\|^{2p}\right] & \leq   C_8 \left(\Delta t\right)^{2p}.
\end{split}
\end{equation}
\begin{equation}
\begin{split}
\label{approx_Yj_2}
\mathbb{E}\left[ \left\| Y^{j}_{t,s} -  Id_n -   \Delta W^j_s  \partial \sigma^j\left(\bar{X}^{j}_{t,\hat{\tau}_{t}}\right) - \frac12 \left(\Delta W^j_s \right)^2 \left(\partial^2 \sigma^j \odot \sigma^j + \left(\partial \sigma^j\right)^2\right)\left( \bar{X}^{j}_{t,\hat{\tau}_{t}}\right) \right\|^{2p}\right] & \leq   C_8 \left(\Delta t\right)^{3p}.
\end{split}
\end{equation}
\end{alem}
\begin{adem}
Let $p \ge 1, t \in [0,T]$ and $s \in [\hat{\tau}_{t},t]$, we recall that $Y^{j}_{t,s}$ is the solution, at time $\Delta W^j_s$, of the following ODE:
\begin{equation}
\left\{
    \begin{array}{ll}
\frac{d}{du} \Xi_{t,\hat{\tau}_{t},u} =  \partial \sigma^j \left( H^{j}_{t,u} \right)   \Xi_{t,\hat{\tau}_{t},u}, u \in \mathbb{R}\\
 \Xi_{t,\hat{\tau}_{t},\hat{\tau}_{t}} = Id_n,
\label{EDOj}
\end{array}
\right.
\end{equation}
where $H^j_{t,u} = h_j\left(\frac{\Delta t}{2},\Delta W^1_t,\ldots,\Delta W^{j-1}_t, u; X^{NV}_{\hat{\tau}_{t}} \right)$.
Since $\partial \sigma^j$ is bounded, there exists a constant $c \in \mathbb{R}_+^*$, independent of $N$, $t$ and $s$, such that
\begin{equation}
\begin{split}
\label{Estim_EDO_Yj}
\left\| Y^{j}_{t,s} \right\|^{2p}   \leq \exp\left(c \left|W^j_s - W^j_{\hat{\tau}_{t}}\right|\right).
\end{split} 
\end{equation}
Therefore, 
\begin{equation}
\begin{split}
\label{Mom_Yj}
\mathbb{E}\left[ \left\| Y^{j}_{t,s}\right\|^{2p} \right]   \leq 2 \exp\left( \frac12 c^2 \left(s - \hat{\tau}_{t}\right)\right) \leq 2\exp\left( \frac12 c^2 T\right).
\end{split} 
\end{equation}
Now, we are able to prove \eqref{Acc_Yj}. Writing $Y^{j}_{t,s}$  in integral form, we have:
\begin{equation*}
\begin{split}
Y^{j}_{t,s} &= Id_n + \int_{0}^{\Delta W^j_s} \displaystyle \partial \sigma^j\left( H^{j}_{t,u}\right) \bar{Y}^{j}_{t,u} du 
\end{split}
\end{equation*}
where $\bar{Y}^{j}_{t,u}$ is the solution at time $u \in \mathbb{R}$ to the ODE \eqref{EDOj}. Applying It\^o's formula, we obtain
\begin{equation}
\begin{split}
\label{_Yj_2}
Y^{j}_{t,s} - Id_n &= \int_{\hat{\tau}_{t}}^{s} \displaystyle \partial \sigma^j\left( \bar{X}^{j}_{t,u}\right) Y^{j}_{t,u} dW^j_u +  \frac12 \int_{\hat{\tau}_{t}}^{s} \displaystyle \left(\partial^2 \sigma^j \odot \sigma^j + \left(\partial \sigma^j\right)^2\right)\left( \bar{X}^{j}_{t,u}\right) Y^{j}_{t,u} du.
\end{split}
\end{equation}
Note that as $\partial \sigma^j \sigma^j$ is Lipschitz continuous, its derivatives given by $\partial^2 \sigma^j \odot \sigma^j + \left(\partial \sigma^j\right)^2$ is bounded. Thus, combining a convexity inequality, the Burkholder-Davis-Gundy inequality and \eqref{Mom_Yj} we deduce \eqref{Acc_Yj}. We prove now \eqref{approx_Yj}. 
Using \eqref{_Yj_2}, we have that

\begin{equation}
\begin{split}
\label{dec_Yj}
Y^{j}_{t,s} -  Id_n -   \Delta W^j_s  \partial \sigma^j\left(\bar{X}^{j}_{t,\hat{\tau}_{t}}\right) &= \int_{\hat{\tau}_{t}}^{s} \displaystyle \left(\partial \sigma^j\left( \bar{X}^{j}_{t,u}\right) Y^{j}_{t,u}  - \partial \sigma^j\left(\bar{X}^{j}_{t,\hat{\tau}_{t}}\right) \right) dW^j_u \\
&+ \frac12  \int_{\hat{\tau}_{t}}^{s} \displaystyle \left(\partial^2 \sigma^j \odot \sigma^j + \left(\partial \sigma^j\right)^2\right)\left( \bar{X}^{j}_{t,u}\right) Y^{j}_{t,u} du \\
&  = \int_{\hat{\tau}_{t}}^{s} \displaystyle \left(\partial \sigma^j\left( \bar{X}^{j}_{t,u}\right)   - \partial \sigma^j\left(\bar{X}^{j}_{t,\hat{\tau}_{t}}\right) \right) Y^{j}_{t,u}  dW^j_u +  \int_{\hat{\tau}_{t}}^{s} \displaystyle  \partial \sigma^j\left(\bar{X}^{j}_{t,\hat{\tau}_{t}}\right) \left(Y^{j}_{t,u} - Id_n \right)  dW^j_u \\
&+   \frac12 \int_{\hat{\tau}_{t}}^{s} \displaystyle \left(\partial^2 \sigma^j \odot \sigma^j + \left(\partial \sigma^j\right)^2\right)\left( \bar{X}^{j}_{t,u}\right) Y^{j}_{t,u} du.
\end{split}
\end{equation}
Since $\partial^2 \sigma^j \odot \sigma^j + \left(\partial \sigma^j\right)^2$ is bounded and $\partial \sigma^j$ is bounded and locally Lipschitz with polynomially growing Lipschitz constant, combining a convexity inequality, the Burkholder-Davis-Gundy inequality, \eqref{rem2_1} from Lemma \ref{Lemme0_A}, \eqref{Acc_Yj}  from Lemma \ref{Lemme4_A} and \eqref{Mom_Yj}, we obtain \eqref{approx_Yj}.  
We can now focus on the approximation \eqref{approx_Yj_2} , with strong order $3/2$, of $Y^j$. We denote by $\Theta^j$ the process such that
\begin{equation}
\Theta^j_{t,s} = Y^{j}_{t,s} -  Id_n -   \Delta W^j_s  \partial \sigma^j\left(\bar{X}^{j}_{t,\hat{\tau}_{t}}\right) - \frac12 \left(\Delta W^j_s \right)^2 \left(\partial^2 \sigma^j \odot \sigma^j + \left(\partial \sigma^j\right)^2\right)\left( \bar{X}^{j}_{t,\hat{\tau}_{t}}\right)
\end{equation}
Writing 
\begin{equation}
\begin{split}
\frac12 \left(\Delta W^j_s \right)^2 \left(\partial^2 \sigma^j \odot \sigma^j + \left(\partial \sigma^j\right)^2\right)\left( \bar{X}^{j}_{t,\hat{\tau}_{t}}\right) &= \int_{\hat{\tau}_{t}}^{s} \displaystyle \Delta W^j_u  \left(\partial^2 \sigma^j \odot \sigma^j + \left(\partial \sigma^j\right)^2\right)\left( \bar{X}^{j}_{t,\hat{\tau}_{t}}\right)   dW^j_u\\
& + \frac12 \int_{\hat{\tau}_{t}}^{s} \displaystyle  \left(\partial^2 \sigma^j \odot \sigma^j + \left(\partial \sigma^j\right)^2\right)\left( \bar{X}^{j}_{t,\hat{\tau}_{t}}\right)   du
\end{split}
\end{equation}
and using \eqref{dec_Yj}, we get
\begin{equation}
\begin{split}
\Theta^j_{t,s} &=    \int_{\hat{\tau}_{t}}^{s} \displaystyle \left(\left(\partial \sigma^j\left( \bar{X}^{j}_{t,u}\right)   - \partial \sigma^j\left(\bar{X}^{j}_{t,\hat{\tau}_{t}}\right) \right) Y^{j}_{t,u}   - \Delta W^j_u  \partial^2 \sigma^j \odot \sigma^j\left(\bar{X}^{j}_{t,\hat{\tau}_{t}}\right)  \right) dW^j_u \\
&+  \int_{\hat{\tau}_{t}}^{s} \displaystyle  \left(\partial \sigma^j\left(\bar{X}^{j}_{t,\hat{\tau}_{t}}\right) \left(Y^{j}_{t,u} - Id_n \right) - \Delta W^j_u \left(\partial \sigma^j\right)^2\left( \bar{X}^{j}_{t,\hat{\tau}_{t}}\right)   \right) dW^j_u \\
&+   \frac12 \int_{\hat{\tau}_{t}}^{s} \displaystyle \left(\left(\partial^2 \sigma^j \odot \sigma^j + \left(\partial \sigma^j\right)^2\right)\left( \bar{X}^{j}_{t,u}\right) Y^{j}_{t,u} - \left(\partial^2 \sigma^j \odot \sigma^j + \left(\partial \sigma^j\right)^2\right)\left( \bar{X}^{j}_{t,\hat{\tau}_{t}}\right)  \right)du.
\end{split}
\end{equation}
Adding and subtracting some appropriate terms, we obtain:
\begin{equation}
\begin{split}
\label{Add_Sub0}
\Theta^j_{t,s} &=    \int_{\hat{\tau}_{t}}^{s} \displaystyle \left(\partial \sigma^j\left( \bar{X}^{j}_{t,u}\right)   - \partial \sigma^j\left(\bar{X}^{j}_{t,\hat{\tau}_{t}}\right)     - \Delta W^j_u  \partial^2 \sigma^j \odot \sigma^j\left(\bar{X}^{j}_{t,\hat{\tau}_{t}}\right)  \right) Y^{j}_{t,u} dW^j_u \\
& + \int_{\hat{\tau}_{t}}^{s} \displaystyle \Delta W^j_u  \partial^2 \sigma^j \odot \sigma^j\left(\bar{X}^{j}_{t,\hat{\tau}_{t}}\right)  \left( Y^{j}_{t,u} - Id_n\right) dW^j_u \\
&+  \int_{\hat{\tau}_{t}}^{s} \displaystyle  \partial \sigma^j\left(\bar{X}^{j}_{t,\hat{\tau}_{t}}\right) \left(Y^{j}_{t,u} - Id_n - \Delta W^j_u \partial \sigma^j\left( \bar{X}^{j}_{t,\hat{\tau}_{t}}\right) \right)    dW^j_u \\
&+   \frac12 \int_{\hat{\tau}_{t}}^{s} \displaystyle  \left( \left(\partial^2 \sigma^j \odot \sigma^j + \left(\partial \sigma^j\right)^2\right)\left( \bar{X}^{j}_{t,u}\right)  - \left(\partial^2 \sigma^j \odot \sigma^j + \left(\partial \sigma^j\right)^2\right)\left( \bar{X}^{j}_{t,\hat{\tau}_{t}}\right)  \right)Y^{j}_{t,u} du\\
& + \frac12 \int_{\hat{\tau}_{t}}^{s} \displaystyle  \left(\partial^2 \sigma^j \odot \sigma^j + \left(\partial \sigma^j\right)^2\right)\left( \bar{X}^{j}_{t,\hat{\tau}_{t}}\right)  \left( Y^{j}_{t,u} - Id_n\right) du.
\end{split}
\end{equation}
On the one hand, by the mean value theorem
\begin{equation}
\label{_Etap_1_Yj}
\partial \sigma^j\left( \bar{X}^j_{t,u}\right) - \partial \sigma^j\left(\bar{X}^j_{t,\hat{\tau}_t}\right) = \partial^2 \sigma^j\left(\bar{X}^j_{t,\hat{\tau}_t}\right) \odot \left( \bar{X}^j_{t,u} - \bar{X}^j_{t,\hat{\tau}_t} \right) + \left(\partial^2 \sigma^j\left(\xi^j_{t,u}\right) -  \partial^2 \sigma^j\left(\bar{X}^j_{t,\hat{\tau}_t}\right) \right) \odot \left( \bar{X}^j_{t,u} - \bar{X}^j_{t,\hat{\tau}_t} \right)
\end{equation}
where $\xi^j_{t,u}$ is a matrix of intermediate points between $\bar{X}^j_{t,u}$ and $\bar{X}^j_{t,\hat{\tau}_t}$. On the other hand, since 
\begin{equation}
\bar{X}^j_{t,u} = \exp\left(\Delta W^j_u \sigma^j \right)\bar{X}^j_{t,\hat{\tau}_t}, 
\end{equation}
applying It\^o's formula we have that:
\begin{equation}
\label{_Etap_2_Yj}
\bar{X}^j_{t,u} - \bar{X}^j_{t,\hat{\tau}_t} = \Delta W^j_u \sigma^j\left(\bar{X}^j_{t,\hat{\tau}_t}\right) + \gamma_{t,u}^j, 
\end{equation}
where 
\begin{equation}
\gamma_{t,u}^j = \int_{\hat{\tau}_{t}}^{u} \left(\sigma^j\left( \bar{X}^j_{t,v}\right) - \sigma^j\left(\bar{X}^j_{t,\hat{\tau}_t}\right) \right)dW^j_v +  \frac12 \int_{\hat{\tau}_{t}}^{u}\partial \sigma^j \sigma^j \left( \bar{X}^j_{t,v}\right) dv.
\end{equation}
Using \eqref{rem2_1} from Lemma \ref{Lemme0_A}, together with the regularity  assumptions on $\sigma^j$ and $\partial \sigma^j \sigma^j$ it is easy to see that $\gamma^j$ a remainder with strong order $1$: $\exists \alpha \in \mathbb{R}_+^*$  independent of $N$, $t$ and $u$, such that
\begin{equation}
\label{gamma}
\mathbb{E}\left[ \left\| \gamma_{t,u}^j  \right\|^{2p} \right]   \leq \alpha \left(\Delta u\right)^{2p}.
\end{equation}
Combining \eqref{Add_Sub0}, \eqref{_Etap_1_Yj} and \eqref{_Etap_2_Yj}, we get 
\begin{equation}
\begin{split}
\Theta^j_{t,s} &=    \int_{\hat{\tau}_{t}}^{s} \displaystyle \left(\partial^2 \sigma^j\left(\bar{X}^j_{t,\hat{\tau}_t}\right) \odot \gamma_{t,u}^j + \left(\partial^2 \sigma^j\left(\xi^j_{t,u}\right) -  \partial^2 \sigma^j\left(\bar{X}^j_{t,\hat{\tau}_t}\right) \right) \odot \left( \bar{X}^j_{t,u} - \bar{X}^j_{t,\hat{\tau}_t} \right)    \right) Y^{j}_{t,u} dW^j_u \\
& + \int_{\hat{\tau}_{t}}^{s} \displaystyle \Delta W^j_u  \partial^2 \sigma^j \odot \sigma^j\left(\bar{X}^{j}_{t,\hat{\tau}_{t}}\right)  \left( Y^{j}_{t,u} - Id_n\right) dW^j_u \\
&+  \int_{\hat{\tau}_{t}}^{s} \displaystyle  \partial \sigma^j\left(\bar{X}^{j}_{t,\hat{\tau}_{t}}\right) \left(Y^{j}_{t,u} - Id_n - \Delta W^j_s \partial \sigma^j\left( \bar{X}^{j}_{t,\hat{\tau}_{t}}\right) \right)    dW^j_u \\
&+   \frac12 \int_{\hat{\tau}_{t}}^{s} \displaystyle  \left( \left(\partial^2 \sigma^j \odot \sigma^j + \left(\partial \sigma^j\right)^2\right)\left( \bar{X}^{j}_{t,u}\right)  - \left(\partial^2 \sigma^j \odot \sigma^j + \left(\partial \sigma^j\right)^2\right)\left( \bar{X}^{j}_{t,\hat{\tau}_{t}}\right)  \right)Y^{j}_{t,u} du\\
& + \frac12 \int_{\hat{\tau}_{t}}^{s} \displaystyle  \left(\partial^2 \sigma^j \odot \sigma^j + \left(\partial \sigma^j\right)^2\right)\left( \bar{X}^{j}_{t,\hat{\tau}_{t}}\right)  \left( Y^{j}_{t,u} - Id_n\right) du.
\end{split}
\end{equation}
Since $\sigma^j$ is Lipschitz continuous, $\partial \sigma^j$ and $\partial^2 \sigma^j$ are locally Lipschitz with polynomially growing Lipschitz constant, it is easy now to see that \eqref{approx_Yj_2} is a straightforward consequence of \eqref{rem2_1} from Lemma \ref{Lemme0_A}, \eqref{Acc_Yj}, \eqref{approx_Yj} and \eqref{gamma}.
\end{adem}

\begin{_adem} \textbf{of Proposition \ref{Prop_App_2}:}\\\\
We recall the result of this proposition:
\begin{equation*}
E = \mathbb{E}\left[ \underset{t\leq T}{\sup} \left\| \displaystyle \int_{0}^{t} \left( Z^0_s - \theta_s^0\right) ds \right\|^{2p}\right] \leq C_2 h^{3p},
\end{equation*}
where, for $t\in [0,T]$, 
\begin{equation*}
Z^0_t =  Y^{d+1}_{t,t} \ldots  Y^{1}_{t,t} \sigma^0\left(\bar{X}_t^0 \right),
\end{equation*}
and 
\begin{equation*}
\begin{split}
\theta_t^0 &=   \sigma^0\left(\bar{X}^0_t\right) + \sum \limits_{j=1}^d  \Delta W^j_t \partial \sigma^j \sigma^0\left(X^{NV}_{\hat{\tau}_t}\right)  + \frac{1}{2} \sum \limits_{j=1}^d \Delta t\left(\partial^2 \sigma^j \odot \sigma^j + \left(\partial \sigma^j\right)^2\right) \sigma^0\left(X^{NV}_{\hat{\tau}_t}\right)  + \frac{1 }{2} \Delta t \partial \sigma^0 \sigma^0\left(X^{NV}_{\hat{\tau}_t}\right)
\end{split}
\end{equation*}
To prove this proposition, we will proceed in two steps. \\\\
\textbf{Step 1: approximation with strong order $3/2$ of $Z^0$.}\\\\
The first step consists in naively computing and approximating the product:   
\begin{equation*}
Z^0_t =  Y^{d+1}_{t,t} \ldots  Y^{1}_{t,t} \sigma^0\left(\bar{X}_t^0 \right).
\end{equation*}                                          
We replace $ Y^{d+1}$ and $Y^{j}$, for $j\in \left\{1,\dots,d\right\}$, by their approximation \eqref{approx_1} from Lemma \ref{Prop_App_1} and \eqref{approx_Yj_2} from Lemma \ref{Lemme4_A}, respectively.
Let
\begin{equation*}
\begin{split}
\Gamma_t &= \left( Id_n +  \frac12 \Delta t  \partial \sigma^0\left(X^{NV}_{t}\right) \right) \left( Id_n +  \Delta W^d_t  \partial \sigma^d\left(\bar{X}^{d}_{t,\hat{\tau}_{t}}\right) + \frac12 \left(\Delta W^d_t \right)^2 \left(\partial^2 \sigma^d \odot \sigma^d + \left(\partial \sigma^d\right)^2\right)\left( \bar{X}^{d}_{t,\hat{\tau}_{t}}\right) \right)\ldots\\
&\left( Id_n +  \Delta W^1_t  \partial \sigma^1\left(\bar{X}^{1}_{t,\hat{\tau}_{t}}\right) + \frac12 \left(\Delta W^1_t \right)^2 \left(\partial^2 \sigma^1 \odot \sigma^1 + \left(\partial \sigma^1\right)^2\right)\left( \bar{X}^{1}_{t,\hat{\tau}_{t}}\right) \right),
\end{split}
\end{equation*} 
and 
\begin{equation*}
\begin{split}
\tilde{\theta}^0_t = \Gamma_t  \sigma^0\left(\bar{X}_t^0 \right).
\end{split}
\end{equation*}
It is clear that $\tilde{\theta}^0$ is an approximation of strong order $3/2$ of $Z^0$: $\exists \tilde{\alpha} \in \mathbb{R}^*_+$ independent of $N$ and $t$ such that:
\begin{equation}
\label{Z0_step1_1} 
\mathbb{E}\left[  \left\|   Z^0_t - \tilde{\theta}_t^0  \right\|^{2p}\right] \leq \tilde{\alpha} h^{3p}.
\end{equation}
Computing $\tilde{\theta}^0$, we easily deduce the following approximation with strong order $3/2$ by dropping higher order terms:
\begin{equation} 
\mathbb{E}\left[  \left\| \hat{\theta}_t^0 - \tilde{\theta}_t^0  \right\|^{2p}\right] \leq \hat{\alpha} h^{3p},
\label{Z0_step1_2}
\end{equation}
where
\begin{equation}
\begin{split}
\hat{\theta}^0_t  &= \sigma^0\left(\bar{X}^0_t\right) + \sum \limits_{j=1}^d  \Delta W^j_t \partial \sigma^j\left(\bar{X}^j_{t,\hat{\tau}_t}\right)\sigma^0\left(\bar{X}^0_t\right)  + \sum \limits_{j=1}^d \sum \limits_{m=1}^{j-1}\Delta W^j_t \Delta W^m_t \partial \sigma^j\left(\bar{X}^j_{t,\hat{\tau}_t}\right) \partial \sigma^m\left(\bar{X}^m_{t,\hat{\tau}_t}\right)\sigma^0\left(\bar{X}^0_t\right) \\
&+  \sum \limits_{j=1}^d \frac{1}{2} \left(\Delta W^j_t\right)^2\left(\partial^2 \sigma^j \odot \sigma^j + \left(\partial \sigma^j\right)^2\right)\left(\bar{X}^j_{t,\hat{\tau}_t}\right) \Bigg) \sigma^0\left(\bar{X}^0_t\right) - \frac{1}{2} \Delta t \partial \sigma^0\left(X^{NV}_t\right)\sigma^0\left(\bar{X}^0_t\right) ,
\end{split}
\end{equation}
and $\hat{\alpha} \in \mathbb{R}^*_+$ is a constant independent of $N$ and $t$.
We will now approximate the term $ \Delta W^j_t \partial \sigma^j\left(\bar{X}^j_{t,\hat{\tau}_t}\right)\sigma^0\left(\bar{X}^0_t\right)$, for $j \in \left\{1,\ldots,d\right\}$, with strong order $3/2$. To do so, it suffices to approximate $\partial \sigma^j\left(\bar{X}^j_{t,\hat{\tau}_t}\right)\sigma^0\left(\bar{X}^0_t\right)$ with strong order $1$.
By the mean value theorem
\begin{equation*}
\begin{split}
\partial \sigma^j\left(\bar{X}^j_{t,\hat{\tau}_t}\right) &= \partial \sigma^j\left(X^{NV}_{\hat{\tau}_t}\right) + \partial^2 \sigma^j\left(X^{NV}_{\hat{\tau}_t}\right) \odot \left(\bar{X}^j_{t,\hat{\tau}_t} - X^{NV}_{\hat{\tau}_t}\right) +  \left( \partial^2 \sigma^j\left(\zeta^{j}_{t}\right)- \partial^2 \sigma^j\left(X^{NV}_{\hat{\tau}_t}\right) \right)\odot \left(\bar{X}^j_{t,\hat{\tau}_t} - X^{NV}_{\hat{\tau}_t}\right),
\end{split}
\end{equation*}
where $\zeta^{j}_{t}$ is a matrix of intermediate points between $\bar{X}^j_{t,\hat{\tau}_t}$ and $X^{NV}_{\hat{\tau}_t}$. On the one hand, since $\forall t \in [0,T]$, $\bar{X}^0_{t} = \bar{X}^{1}_{t,\hat{\tau}_t}$ and $\bar{X}^{m+1}_{t,\hat{\tau}_t} = \bar{X}^{m}_{t,t}$, using telescopic summation, we have
\begin{equation}
\begin{split}
\label{step___1}
\partial \sigma^j\left(\bar{X}^j_{t,\hat{\tau}_t}\right) &=  \partial \sigma^j\left(X^{NV}_{\hat{\tau}_t}\right)  + \sum \limits_{m=1}^{j-1}\left( \partial^2 \sigma^j\left(X^{NV}_{\hat{\tau}_t}\right) \odot \left(\bar{X}^{m+1}_{t,\hat{\tau}_t} - \bar{X}^{m}_{t,\hat{\tau}_t}\right) \right)\\
&+ \partial^2 \sigma^j\left(X^{NV}_{\hat{\tau}_t}\right) \odot \left(\bar{X}^0_{t} - X^{NV}_{\hat{\tau}_t}\right) + \left( \partial^2 \sigma^j\left(\zeta^{j}_{t}\right)- \partial^2 \sigma^j\left(X^{NV}_{\hat{\tau}_t}\right) \right)\odot \left(\bar{X}^j_{t,\hat{\tau}_t} - X^{NV}_{\hat{\tau}_t}\right)\\
& = \partial \sigma^j\left(X^{NV}_{\hat{\tau}_t}\right)  + \sum \limits_{m=1}^{j-1}\left( \partial^2 \sigma^j\left(X^{NV}_{\hat{\tau}_t}\right) \odot \left(\bar{X}^{m}_{t,t} - \bar{X}^{m}_{t,\hat{\tau}_t}\right) \right)\\
&+ \partial^2 \sigma^j\left(X^{NV}_{\hat{\tau}_t}\right) \odot \left(\bar{X}^0_{t} - X^{NV}_{\hat{\tau}_t}\right) + \left( \partial^2 \sigma^j\left(\zeta^{j}_{t}\right)- \partial^2 \sigma^j\left(X^{NV}_{\hat{\tau}_t}\right) \right)\odot \left(\bar{X}^j_{t,\hat{\tau}_t} - X^{NV}_{\hat{\tau}_t}\right).
\end{split}
\end{equation}
On the other hand, we recall that from \eqref{_Etap_2_Yj} we have 
\begin{equation}
\bar{X}^m_{t,t} - \bar{X}^m_{t,\hat{\tau}_t} = \Delta W^m_t \sigma^m\left(\bar{X}^m_{t,\hat{\tau}_t}\right) + \gamma_{t,t}^m, 
\end{equation}
where $\gamma^m$ is a remainder with strong order $1$:
\begin{equation}
\gamma_{t,t}^m = \int_{\hat{\tau}_{t}}^{t} \left(\sigma^m\left( \bar{X}^m_{t,v}\right) - \sigma^m\left(\bar{X}^m_{t,\hat{\tau}_t}\right) \right)dW^m_v +  \frac12 \int_{\hat{\tau}_{t}}^{t}\partial \sigma^m \sigma^m \left( \bar{X}^m_{t,v}\right) dv.
\end{equation}
This leads us to the following decomposition 
\begin{equation}
\begin{split}
\label{step___1__}
\partial \sigma^j\left(\bar{X}^j_{t,\hat{\tau}_t}\right) &=   \partial \sigma^j\left(X^{NV}_{\hat{\tau}_t}\right)  + \sum \limits_{m=1}^{j-1} \Delta W^m_t  \partial^2 \sigma^j\left(X^{NV}_{\hat{\tau}_t}\right) \odot \sigma^m\left(\bar{X}^m_{t,\hat{\tau}_t}\right) + \Gamma^j_t,
\end{split}
\end{equation}
where
\begin{equation}
\begin{split}
\Gamma^j_t &= \sum \limits_{m=1}^{j-1} \Delta W^m_t  \partial^2 \sigma^j\left(X^{NV}_{\hat{\tau}_t}\right) \odot \gamma_{t,t}^m + \partial^2 \sigma^j\left(X^{NV}_{\hat{\tau}_t}\right) \odot \left(\bar{X}^0_{t} - X^{NV}_{\hat{\tau}_t}\right) \\
&+ \left( \partial^2 \sigma^j\left(\zeta^{j}_{t}\right)- \partial^2 \sigma^j\left(X^{NV}_{\hat{\tau}_t}\right) \right)\odot \left(\bar{X}^j_{t,\hat{\tau}_t} - X^{NV}_{\hat{\tau}_t}\right).
\end{split}
\end{equation}
Using \eqref{Acc_0} from Lemma \ref{Lemme0} and \eqref{Acc_j_tau} from Lemma \ref{Lemme0_A}, it is easy to see that $\Gamma^j$ is a remainder with strong order $1$. Then it follows that $\check{\theta}^0$, which is given by
\begin{equation}
\begin{split}
\check{\theta}^0_t &= \sigma^0\left(\bar{X}^0_t\right) + \sum \limits_{j=1}^d  \Delta W^j_t \partial \sigma^j\left(X^{NV}_{\hat{\tau}_t}\right)\sigma^0\left(\bar{X}^0_t\right)  + \sum \limits_{j=1}^d \sum \limits_{m=1}^{j-1}\Delta W^j_t \Delta W^m_t \partial \sigma^j\left(\bar{X}^j_{t,\hat{\tau}_t}\right) \partial \sigma^m\left(\bar{X}^m_{t,\hat{\tau}_t}\right)\sigma^0\left(\bar{X}^0_t\right) \\
&+ \sum \limits_{j=1}^d \sum \limits_{m=1}^{j-1}\Delta W^j_t \Delta W^m_t \partial^2 \sigma^j \left(X^{NV}_{\hat{\tau}_t}\right) \odot \sigma^m\left(\bar{X}^m_{t,\hat{\tau}_t}\right)\sigma^0\left(\bar{X}^0_t\right) \\
&+  \sum \limits_{j=1}^d \frac{1}{2} \left(\Delta W^j_t\right)^2\left(\partial^2 \sigma^j \odot \sigma^j + \left(\partial \sigma^j\right)^2\right)\left(\bar{X}^j_{t,\hat{\tau}_t}\right) \Bigg) \sigma^0\left(\bar{X}^0_t\right) - \frac{1}{2} \Delta t \partial \sigma^0\left(X^{NV}_t\right)\sigma^0\left(\bar{X}^0_t\right) ,
\end{split}
\end{equation}
is an approximation of strong order $3/2$ of $\hat{\theta}^0$:$\exists \check{\alpha} \in \mathbb{R}^*_+$ independent of $N$ and $t$ such that:
\begin{equation} 
\mathbb{E}\left[  \left\|   \check{\theta}_t^0 - \hat{\theta}_t^0  \right\|^{2p}\right] \leq \check{\alpha} h^{3p}.
\label{Z0_step1_3}
\end{equation}
Then, we can replace by $X^{NV}_{\hat{\tau}_t}$ the argument of the functions in the expression of $\check{\theta}^0$  when Lemma \ref{Lemme0} and Lemma \ref{Lemme0_A} ensure that the strong order $3/2$ is preserved. This leads us to the following approximation, given by $\bar{\theta}^0$,
\begin{equation}
\begin{split}
\label{bar_theta_0}
\bar{\theta}_t^0 &=   \sigma^0\left(\bar{X}^0_t\right) + \sum \limits_{j=1}^d  \Delta W^j_t \partial \sigma^j \sigma^0\left(X^{NV}_{\hat{\tau}_t}\right)  + \sum \limits_{j=1}^d \sum \limits_{m=1}^{j-1}\Delta W^j_t \Delta W^m_t \left(\partial^2 \sigma^j \odot \sigma^m + \partial \sigma^j \partial \sigma^m  \right) \sigma^0\left(X^{NV}_{\hat{\tau}_t}\right)\\
&+ \sum \limits_{j=1}^d \frac{1}{2} \left(\Delta W^j_t\right)^2\left(\partial^2 \sigma^j \odot \sigma^j + \left(\partial \sigma^j\right)^2\right) \sigma^0\left(X^{NV}_{\hat{\tau}_t}\right)  + \frac{1 }{2} \Delta t \partial \sigma^0 \sigma^0\left(X^{NV}_{\hat{\tau}_t}\right), 
\end{split}
\end{equation}
and $ \exists \bar{\alpha} \in \mathbb{R}^*_+$ independent of $N$ and $t$, such that
\begin{equation} 
\label{Z0_step1_4}
\mathbb{E}\left[  \left\|   \bar{\theta}_t^0 - \check{\theta}_t^0  \right\|^{2p}\right] \leq \bar{\alpha} h^{3p}.
\end{equation}
Adding and subtracting $ \bar{\theta}^0$, and using a convexity inequality, we get:
\begin{equation}
\begin{split}
E = \mathbb{E}\left[ \underset{t\leq T}{\sup} \left\| \displaystyle \int_{0}^{t} \left( Z^0_s - \theta_s^0\right) ds \right\|^{2p}\right] \leq 1 \vee T^{2p-1} \left( \displaystyle \int_{0}^{T} \mathbb{E}\left[  \left\|  Z^0_t - \bar{\theta}_t^0  \right\|^{2p}\right]dt +  \mathbb{E}\left[ \underset{t\leq T}{\sup} \left\| \displaystyle \int_{0}^{t} \left(\bar{\theta}_s^0 - \theta_s^0\right) ds \right\|^{2p}\right] \right).
\end{split}
\end{equation}
Combining a convexity inequality, \eqref{Z0_step1_1}, \eqref{Z0_step1_2}, \eqref{Z0_step1_3}, and \eqref{Z0_step1_4}, we obtain an estimation of the first expectation in the right-hand side of the above inequality:
\begin{equation}
\begin{split}
\mathbb{E}\left[  \left\|  Z^0_t - \bar{\theta}_t^0  \right\|^{2p}\right] \leq 4^{2p-1} \left(\tilde{\alpha} + \hat{\alpha} + \check{\alpha} + \bar{\alpha} \right)h^{3p}
\end{split}
\end{equation}
To achieve our goal, it remains to estimate:
\begin{equation}
\mathbb{E}\left[ \underset{t\leq T}{\sup} \left\| \displaystyle \int_{0}^{t} \left(\bar{\theta}_s^0 - \theta_s^0\right) ds \right\|^{2p}\right]. 
\end{equation} 
This is the aim of the following step.\\\\
\textbf{Step 2: the integration by parts formula.}\\\\
Let $t \in [0,T]$, subtracting \eqref{sub_theta} from  \eqref{bar_theta_0}, we have that
\begin{equation}
\begin{split}
\displaystyle \int_{0}^{t} \left(\bar{\theta}_s^0 - \theta_s^0 \right)ds &= \sum \limits_{j=1}^d \sum \limits_{m=1}^{j-1} \displaystyle \int_{0}^{t}  \Delta W^j_s \Delta W^m_s  \left(\partial^2 \sigma^j \odot \sigma^m + \partial \sigma^j \partial \sigma^m  \right)  \sigma^0\left(X^{NV}_{\hat{\tau}_s}\right) ds\\
& + \frac12 \displaystyle \int_{0}^{t} \left(  \left(\Delta W^j_s\right)^2 - \Delta s \right)  \left(\partial^2 \sigma^j \odot \sigma^j + \left(\partial \sigma^j  \right)^2 \right)  \sigma^0\left(X^{NV}_{\hat{\tau}_s}\right) ds.
\end{split}
\end{equation}
To lighten up the notation, we denote $ F_{j,m} = \left(\partial^2 \sigma^j \odot \sigma^m + \partial \sigma^j \partial \sigma^m  \right) \sigma^0$ for $ j,m \in \left\{1,\ldots,d\right\}, m \leq j$.
Using the integration by parts formula, we have:
\begin{equation}
\displaystyle \int_{0}^{t}  \Delta W^j_s \Delta W^m_s  F_{j,m}\left(X^{NV}_{\hat{\tau}_s}\right) ds = \displaystyle \int_{0}^{t}  \left(t \wedge  \check{\tau}_s - s\right) \Delta W^j_s  F_{j,m}\left(X^{NV}_{\hat{\tau}_s}\right) dW^m_s  + \displaystyle \int_{0}^{t}  \left(t \wedge  \check{\tau}_s - s\right) \Delta W^m_s  F_{j,m}\left(X^{NV}_{\hat{\tau}_s}\right) dW^j_s.
\end{equation}
Taking the expectation of the supremum and using a convexity inequality, the Burkholder-Davis-Gundy inequality, we easily get a positive constant $\alpha_1$, independent of $N$ such that
\begin{equation}
\mathbb{E}\left[ \underset{t\leq T}{\sup} \left\| \displaystyle \int_{0}^{t}  \Delta W^j_s \Delta W^m_s  F_{j,m}\left(X^{NV}_{\hat{\tau}_s}\right) ds  \right\|^{2p}\right] \leq \alpha_1  \displaystyle \int_{0}^{T} \left(\check{\tau}_s - s\right)^{2p} \mathbb{E}\left[\left\| \Delta W^m_s  F_{j,m}\left(X^{NV}_{\hat{\tau}_s}\right)  \right\|^{2p}\right] ds .
\end{equation}
Then, by independence
\begin{equation}
\mathbb{E}\left[ \underset{t\leq T}{\sup} \left\| \displaystyle \int_{0}^{t}  \Delta W^j_s \Delta W^m_s  F_{j,m}\left(X^{NV}_{\hat{\tau}_s}\right) ds  \right\|^{2p}\right] \leq \alpha_1  \displaystyle \int_{0}^{T} \left(\check{\tau}_s - s\right)^{2p} \mathbb{E}\left[\left| \Delta W^m_s  \right|^{2p}\right] \mathbb{E}\left[\left\| F_{j,m}\left(X^{NV}_{\hat{\tau}_s}\right)  \right\|^{2p}\right] ds. 
\end{equation}
Since $F_{j,m}\in LIP_{loc}^{pgc}\left(\mathbb{R}^n \right)$, we get $\alpha_2\in \mathbb{R}_+$ and $q \ge 1$ independent of $N$ such that
\begin{equation}
\mathbb{E}\left[ \underset{t\leq T}{\sup} \left\| \displaystyle \int_{0}^{t}  \Delta W^j_s \Delta W^m_s  F_{j,m}\left(X^{NV}_{\hat{\tau}_s}\right) ds  \right\|^{2p}\right] \leq \alpha_2 \mathbb{E}\left[\left| G \right|^{2p}\right] \left(\displaystyle \int_{0}^{T} \mathbb{E}\left[\left\| \left(X^{NV}_{\hat{\tau}_s}\right)  \right\|^{2q}\right] ds\right)  h^{3p}
\end{equation}
where $G$ is a standard normal variable. We conclude using \eqref{NV_moment} from Lemma \ref{Lemme0}. 
To estimate, 
\begin{equation}
\displaystyle \int_{0}^{t} \left(  \left(\Delta W^j_s\right)^2 - \Delta s \right)  F_{j,j}\left(X^{NV}_{\hat{\tau}_s}\right) ds,
\end{equation}
we use exactly the same arguments since the integration by parts formula gives us
\begin{equation}
\displaystyle \int_{0}^{t} \left(  \left(\Delta W^j_s\right)^2 - \Delta s \right)  F_{j,j}\left(X^{NV}_{\hat{\tau}_s}\right) ds = 2 \displaystyle \int_{0}^{t}  \left(t \wedge  \check{\tau}_s - s\right) \Delta W^j_s  F_{j,j}\left(X^{NV}_{\hat{\tau}_s}\right) dW^j_s.  
\end{equation}
This completes the proof.
\end{_adem}

In the following, we give a detailed proof of Proposition \ref{Prop_Q}, which is a consequence of Propositions \ref{Prop_App_1}, \ref{Prop_App_2}.

\begin{_adem} \textbf{of Proposition \ref{Prop_Q}:}\\\\
Combining a convexity inequality and the Burkholder-Davis-Gundy, we get a constant $\alpha_0$ independent of $N$ such that
\begin{equation}
\begin{split}
\label{Conve}
\mathbb{E}\left[ \underset{t\leq T}{\sup} \left\| \bar{Q}^{1,N}_t  \right\|^{2p}\right] & \leq  \alpha_0 \left( E + \sum \limits_{j=1}^d \left( E^j_1 + E^j_2 + E^j_3 \right) \right),
\end{split}
\end{equation}
where
\begin{equation}
\begin{split}
E = \mathbb{E}\left[ \underset{t\leq T}{\sup} \left\| \displaystyle \int_{0}^{t} \left( Z^0_s - \theta_s^0\right) ds \right\|^{2p}\right],
\end{split}
\end{equation}
for $j \in \left\{1,\ldots,d\right\}$,
\begin{equation}
\begin{split}
 E^j_1 = \displaystyle \int_{0}^{T} \mathbb{E}\left[ \left\| Y^{d+1}_{t,t} \sigma^j\left(\bar{X}_{t,t}\right) - \sigma^j\left(X^{NV}_t\right) - \frac12 \Delta t \left[ \sigma^j, \sigma^0\right ]\left(X^{NV}_t\right) \right\|^{2p}\right] dt,
\end{split}
\end{equation}
\begin{equation}
\begin{split}
 E^j_2 = \displaystyle \int_{0}^{T} \mathbb{E}\left[ \left\| Y^{d+1}_{t,t} \partial \sigma^j \sigma^j\left(\bar{X}_{t,t}\right) - \partial\sigma^j \sigma^j\left(X^{NV}_t\right)  - \frac12 \Delta t \left( \partial \sigma^0 \partial \sigma^j  \sigma^j - \left(\partial^2 \sigma^j \odot \sigma^j + \left(\partial \sigma^j\right)^2\right) \sigma^0 \right)\left(X^{NV}_{\hat{\tau}_t}\right) \right\|^{2p}\right] dt,
\end{split}
\end{equation}
and
\begin{equation}
\begin{split}
 E^j_3 = \displaystyle \int_{0}^{T} \mathbb{E}\left[ \left\| Z^{d+1,j}_t - \frac12 \Delta t \left( \partial^2\sigma^0 \odot \sigma^j \right)\sigma^j\left(X^{NV}_{\hat{\tau}_t}\right)  \right\|^{2p}\right] dt.
\end{split}
\end{equation}
To prove our claim, it suffices to estimate, with at least strong order $3/2$, each expectations in the right-hand side of \eqref{Conve}.\\\\
\textbf{An estimation with strong order $3/2$ of $E$} is given by Proposition \ref{Prop_App_2}.\\\\
\textbf{Estimation with strong order $2$ of $E_1^j$, for $j \in \left\{1,\ldots,d\right\}$.}\\\\
To estimate $E_1^j$, we estimate its integrand, and we denote 
\begin{equation}
 \epsilon_1^j\left(t\right) = \mathbb{E}\left[ \left\| \lambda_t^j \right\|^{2p}\right],  
\end{equation}
where 
\begin{equation}
\lambda_t^j =  Y^{d+1}_{t,t} \sigma^j\left(\bar{X}_{t,t}\right) - \sigma^j\left(X^{NV}_t\right) - \frac12 \Delta t \left[ \sigma^j, \sigma^0\right ]\left(X^{NV}_t\right).
\end{equation}
Since
\begin{equation}
\begin{split}
\lambda_t^j &= \left(Y^{d+1}_{t,t}  - Id_n - \frac{1}{2} \Delta t  \partial \sigma^0\left(X^{NV}_t\right) \right) \left( \sigma^j\left(\bar{X}_{t,t}\right) - \sigma^j\left(X^{NV}_{t}\right) + \frac{1}{2} \Delta t \partial \sigma^j \sigma^0 \left(X^{NV}_{t} \right) \right)\\
& + \left(Y^{d+1}_{t,t}  - Id_n - \frac{1}{2} \Delta t  \partial \sigma^0\left(X^{NV}_t\right) \right) \left(  \sigma^j\left(X^{NV}_{t}\right) - \frac{1}{2} \Delta t \partial \sigma^j \sigma^0 \left(X^{NV}_{t} \right) \right) \\
&+ \left( Id_n + \frac{1}{2} \Delta t  \partial \sigma^0\left(X^{NV}_t\right) \right) \left( \sigma^j\left(\bar{X}_{t,t}\right) - \sigma^j\left(X^{NV}_{t}\right) + \frac{1}{2} \Delta t \partial \sigma^j \sigma^0 \left(X^{NV}_{t} \right) \right)\\
&- \frac14 \left( \Delta t\right)^2 \partial \sigma^0 \partial \sigma^j \sigma^0\left(X^{NV}_t\right), 
\end{split}
\end{equation}
combining a convexity inequality and the Cauchy-Schwarz inequality, we obtain
\begin{equation}
\begin{split}
 \epsilon_1^j\left(t\right) & \leq  4^{2p-1}\Bigg( \left(\mathbb{E}\left[ \left\| Y^{d+1}_{t,t}  - Id_n - \frac{1}{2} \Delta t  \partial \sigma^0\left(X^{NV}_t\right) \right\|^{4p}\right] \mathbb{E}\left[ \left\|\sigma^j\left(\bar{X}_{t,t}\right) - \sigma^j\left(X^{NV}_{t}\right) + \frac{1}{2} \Delta t \partial \sigma^j \sigma^0 \left(X^{NV}_{t} \right) \right\|^{4p}\right]\right)^{\frac12}\\
&+ \left(\mathbb{E}\left[ \left\| Y^{d+1}_{t,t}  - Id_n - \frac{1}{2} \Delta t  \partial \sigma^0\left(X^{NV}_t\right) \right\|^{4p}\right] \mathbb{E}\left[ \left\| \sigma^j\left(X^{NV}_{t}\right) - \frac{1}{2} \Delta t \partial \sigma^j \sigma^0 \left(X^{NV}_{t} \right) \right\|^{4p}\right]\right)^{\frac12}\\
&+ \left(\mathbb{E}\left[ \left\|  Id_n + \frac{1}{2} \Delta t  \partial \sigma^0\left(X^{NV}_t\right) \right\|^{4p}\right] \mathbb{E}\left[ \left\|\sigma^j\left(\bar{X}_{t,t}\right) - \sigma^j\left(X^{NV}_{t}\right) + \frac{1}{2} \Delta t \partial \sigma^j \sigma^0 \left(X^{NV}_{t} \right) \right\|^{4p}\right]\right)^{\frac12}\\
& + \mathbb{E}\left[ \left\| \frac14 \left( \Delta t\right)^2 \partial \sigma^0 \partial \sigma^j \sigma^0\left(X^{NV}_t\right) \right\|^{2p}\right]\Bigg).
\end{split}
\end{equation}
The first two expectations in the right-hand side of the previous inequality are estimated using \eqref{approx_1} and \eqref{approx_2} from Proposition \ref{Prop_App_2}, respectively. Since $\sigma^0$ and $\sigma^j$ have bounded first order derivatives, applying and \eqref{NV_moment} from Lemma \ref{Lemme0}, we obtain a constant $\beta_j \in \mathbb{R}_+^*$ independent of $N$ such that:  
\begin{equation}
\mathbb{E}\left[ \left\| \sigma^j\left(X^{NV}_{t}\right) - \frac{1}{2} \Delta t \partial \sigma^j \sigma^0 \left(X^{NV}_{t} \right) \right\|^{4p}\right] \leq \beta_j,
\end{equation}
\begin{equation}
\label{estim_preuve}
\mathbb{E}\left[ \left\|  Id_n + \frac{1}{2} \Delta t  \partial \sigma^0\left(X^{NV}_t\right) \right\|^{4p}\right]\leq \beta_j,
\end{equation}
and
\begin{equation}
\mathbb{E}\left[ \left\| \frac14 \left( \Delta t\right)^2 \partial \sigma^0 \partial \sigma^j \sigma^0\left(X^{NV}_t\right) \right\|^{2p}\right] \leq \beta_j h^{4p}.
\end{equation}
Then it follows that
\begin{equation}
\epsilon_1^j\left(t\right) \leq  4^{2p-1} \left(C_1 T^{4p} + 2 \sqrt{C_1 \beta_j } + \beta_j  \right)  h^{4p},
\end{equation}
and then
\begin{equation}
E_1^j  \leq  4^{2p-1} T \left(C_1 T^{4p} + 2 \sqrt{C_1 \beta_j } + \beta_j  \right)  h^{4p}.
\end{equation}

\textbf{Estimation with strong order $3/2$ of $E_2^j$, for $j \in \left\{1,\ldots,d\right\}$.}\\\\
As previously, we denote
\begin{equation}
 \epsilon_2^j\left(t\right) = \mathbb{E}\left[ \left\| \mu_t^j \right\|^{2p}\right],  
\end{equation}
where 
\begin{equation}
\mu_t^j =  Y^{d+1}_{t,t} \partial \sigma^j \sigma^j\left(\bar{X}_{t,t}\right) - \partial\sigma^j \sigma^j\left(X^{NV}_t\right)  - \frac12 \Delta t \left( \partial \sigma^0 \partial \sigma^j \sigma^j - \left(\partial^2 \sigma^j \odot \sigma^j + \left(\partial \sigma^j\right)^2\right) \sigma^0 \right)\left(X^{NV}_{\hat{\tau}_t}\right).
\end{equation}
Since
\begin{equation}
\begin{split}
\mu_t^j &= \left(Y^{d+1}_{t,t}  - Id_n - \frac{1}{2} \Delta t  \partial \sigma^0\left(X^{NV}_t\right) \right) \left( \partial\sigma^j\sigma^j\left(\bar{X}_{t,t}\right) - \partial\sigma^j\sigma^j\left(X^{NV}_{t}\right) + \frac{1}{2} \Delta t \left(\partial^2\sigma^j \odot \sigma^j  + \left(\partial \sigma^j\right)^2 \right) \sigma^0 \left(X^{NV}_{\hat{\tau}_t} \right) \right)\\
& + \left(Y^{d+1}_{t,t}  - Id_n - \frac{1}{2} \Delta t  \partial \sigma^0\left(X^{NV}_t\right) \right) \left(  \partial\sigma^j\sigma^j\left(X^{NV}_{t}\right) - \frac{1}{2} \Delta t \left(\partial^2\sigma^j \odot \sigma^j  + \left(\partial \sigma^j\right)^2 \right) \sigma^0 \left(X^{NV}_{\hat{\tau}_t} \right) \right) \\
&+ \left( Id_n + \frac{1}{2} \Delta t  \partial \sigma^0\left(X^{NV}_t\right) \right) \left( \partial\sigma^j\sigma^j\left(\bar{X}_{t,t}\right) - \partial\sigma^j\sigma^j\left(X^{NV}_{t}\right) + \frac{1}{2} \Delta t \left(\partial^2\sigma^j \odot \sigma^j  + \left(\partial \sigma^j\right)^2 \right) \sigma^0 \left(X^{NV}_{\hat{\tau}_t} \right) \right)\\
&- \frac14 \left( \Delta t\right)^2 \partial \sigma^0 \left(X^{NV}_t\right)\left(\partial^2\sigma^j \odot \sigma^j  + \left(\partial \sigma^j\right)^2 \right) \sigma^0 \left(X^{NV}_{\hat{\tau}_t} \right) + \frac12 \Delta t \left( \partial \sigma^0 \partial \sigma^j \sigma^j\left(X^{NV}_{t} \right) - \partial \sigma^0 \partial \sigma^j \sigma^j\left(X^{NV}_{\hat{\tau}_t} \right) \right) , 
\end{split}
\end{equation}
combining a convexity inequality, the Cauchy-Schwarz inequality, \eqref{approx_1} and \eqref{approx_3} from Proposition \ref{Prop_App_2}, and \eqref{estim_preuve}, we obtain
\begin{equation}
\begin{split}
\label{eps2_preuve}
 \epsilon_2^j\left(t\right) & \leq  5^{2p-1}\Bigg( C_1 h^{7p} + \sqrt{C_1} h^{4p} \left(\mathbb{E}\left[ \left\| \partial\sigma^j\sigma^j\left(X^{NV}_{t}\right) - \frac{1}{2} \Delta t \left(\partial^2\sigma^j \odot \sigma^j  + \left(\partial \sigma^j\right)^2 \right) \sigma^0 \left(X^{NV}_{\hat{\tau}_t} \right) \right\|^{4p}\right] \right)^{\frac12}\\
&+ \sqrt{C_1 \beta_j} h^{3p}  + \mathbb{E}\left[ \left\| \frac14 \left( \Delta t\right)^2 \partial \sigma^0 \left(X^{NV}_t\right)\left(\partial^2\sigma^j \odot \sigma^j  + \left(\partial \sigma^j\right)^2 \right) \sigma^0 \left(X^{NV}_{\hat{\tau}_t} \right) \right\|^{2p}\right]\\
& + \mathbb{E}\left[ \left\| \frac12 \Delta t \left( \partial \sigma^0 \partial \sigma^j \sigma^j\left(X^{NV}_{t} \right) - \partial \sigma^0 \partial \sigma^j \sigma^j\left(X^{NV}_{\hat{\tau}_t} \right) \right)\right\|^{2p}\right]   \Bigg).
\end{split}
\end{equation}
Since $\sigma^0$ and $\sigma^j$ have bounded first order derivatives and $\partial^2 \sigma^j$ is locally Lipschitz with polynomially growing Lipschitz constant, applying and \eqref{NV_moment} from Lemma \ref{Lemme0}, we easily get a constant $\gamma_j \in \mathbb{R}_+^*$ independent of $N$ such that:
\begin{equation}
\mathbb{E}\left[ \left\| \partial\sigma^j\sigma^j\left(X^{NV}_{t}\right) - \frac{1}{2} \Delta t \left(\partial^2\sigma^j \odot \sigma^j  + \left(\partial \sigma^j\right)^2 \right) \sigma^0 \left(X^{NV}_{\hat{\tau}_t} \right) \right\|^{4p}\right] \leq \gamma_j,
\end{equation}
and 
\begin{equation}
\mathbb{E}\left[ \left\| \frac14 \left( \Delta t\right)^2 \partial \sigma^0 \left(X^{NV}_t\right)\left(\partial^2\sigma^j \odot \sigma^j  + \left(\partial \sigma^j\right)^2 \right) \sigma^0 \left(X^{NV}_{\hat{\tau}_t} \right) \right\|^{2p}\right] \leq \gamma_jh^{4p}.
\end{equation}
It remains to estimate the last expectation in the right-hand side of \eqref{eps2_preuve}. As the function $F^j := \displaystyle \frac12 \partial \sigma^0 \partial \sigma^j \sigma^j$ is locally Lipschitz with polynomially growing Lipschitz constant, we easily get a constant $\delta_j \in \mathbb{R}_+^*$ independent of $N$ such that:  
\begin{equation}
\begin{split}  
\mathbb{E}\left[ \left\| \Delta t \left( F^j\left(X^{NV}_t\right) -  F^j\left(X^{NV}_{\hat{\tau}_t} \right)\right) \right\|^{2p}\right]  \leq \delta_j h^{3p}. 
\end{split}
\end{equation}
Then it follows that
\begin{equation}
\epsilon_2^j\left(t\right) \leq  5^{2p-1} \left(C_1 T^{4p} + \sqrt{C_1} T \gamma_j + \sqrt{C_1 \beta_j } + \gamma_j T^p + \delta_j  \right)  h^{3p},
\end{equation}
and we conclude that 
\begin{equation}
E_2^j  \leq  5^{2p-1} T  \left(C_1 T^{4p} + \sqrt{C_1} T \gamma_j + \sqrt{C_1 \beta_j } + \gamma_j T^p + \delta_j  \right)  h^{3p}.
\end{equation}
\textbf{An estimation with strong order $3/2$ of $E_3^j$, for $j \in \left\{1,\ldots,d\right\}$} is given by \eqref{approx_4} from Proposition \ref{Prop_App_1} and this concludes the proof. 
\end{_adem}

\end{document}